\documentclass[12pt]{article}
\usepackage{amsfonts}
\usepackage{graphics}
\newtheorem{THEOREM}{Theorem}[section]
\newtheorem{LEMMA}{Lemma}[section]

\newcommand{\qed}{\ \rule[-1pt]{4pt}{8pt} 

                                         \vspace{2ex} }
\newenvironment{PROOF}{

                       \noindent{\bf Proof}.}{\qed}
\newcounter{labelflag} 

\newcommand{\Label}[1]{
                       \ifnum\thelabelflag=1 
                          \ifmmode  
                             \makebox[0in][l]{\qquad\fbox{\rm#1}}
                          \else
                             \marginpar{\vspace{0.7\baselineskip}
                                        \hspace{-1.1\textwidth}
                                        \fbox{\rm#1}}
                          \fi 
                       \fi
                       \label{#1} 
                      }
\newcommand{\be}{\begin{equation}}
\newcommand{\ee}{\end{equation}}
\newcommand{\eps}{\varepsilon}

\newcommand{\Fe}{{\cal F}_\varepsilon}
\newcommand{\Me}{{\cal M}_\varepsilon}
\newcommand{\Fo}{{\cal F}_0}
\newcommand{\Mo}{{\cal M}_0}

\newcommand{\bfL}{\Lambda}

\newcommand{\Af}{A}
\newcommand{\Bf}{B}
\newcommand{\Uf}{U}
\newcommand{\Lf}{L}
\newcommand{\Kf}{\mathcal{K}_\varepsilon}

\begin{document}

\begin{center}
{\large\textbf{FAST AND SLOW DYNAMICS FOR THE}}\\[1ex]
{\large\textbf{COMPUTATIONAL SINGULAR PERTURBATION}}\\[1ex]
{\large\textbf{METHOD}}\\[1ex]

\vspace{2ex}
Antonios~Zagaris,$^1$
Hans G.~Kaper,$^2$
Tasso J.~Kaper$^1$

$^1$~Department of Mathematics and Center for BioDynamics\\
Boston University,
Boston, Massachusetts, USA \\[1ex]
$^2$~Mathematics and Computer Science Division,\\
Argonne National Laboratory,
Argonne, Illinois, USA \\[1ex]
\end{center}

\begin{small}
\noindent\textbf{Abstract.}
The Computational Singular Perturbation (CSP) method
of Lam and Goussis
is an iterative method to
reduce the dimensionality of
systems of 
ordinary differential equations with 
multiple time scales. 
In [J. Nonlin. Sci., to appear], 
the authors showed that 
each iteration of
the CSP algorithm 
improves the approximation of 
the slow manifold by 
one order. 
In this paper, 
it is shown that
the CSP method
simultaneously approximates 
the tangent spaces to 
the fast fibers 
along which solutions relax to 
the slow manifold. 
Again, each iteration adds 
one order of accuracy.
In some studies,
the output of
the CSP algorithm is 
postprocessed by
linearly projecting initial data onto
the slow manifold along
these approximate tangent spaces.
These projections,
in turn,
also become successively more accurate.



\end{small}

\section{Introduction \label{s-intro}}
\setcounter{equation}{0}
The Computational Singular Perturbation (CSP) method
is one of several so-called reduction methods
developed in chemistry to systematically decrease
the size and complexity of systems of chemical kinetics
equations.
The method was first proposed by
Lam and Goussis~\cite{GL-1992,
L-1993, LG-1988, LG-1991, LG-1994}
and is widely used, for example,
in combustion modeling~\cite{HG-1999, LJL-2001, MDMG-1999, MG-2001, VG-2001, VNG-2003}.

The CSP method is generally applicable
to systems of nonlinear ordinary
differential equations (ODEs)
with simultaneous fast and slow dynamics
where the long-term dynamics evolve on
a low-dimensional slow manifold
in the phase space.
The method is essentially an algorithm
to find successive approximations
to the slow manifold and
match the initial conditions
to the dynamics
on the slow manifold.

In a previous paper~\cite{ZKK-2003},
we focused on the slow manifold
and the accuracy of the CSP approximation
for fast--slow systems of ODEs.
In such systems,
the ratio of the characteristic fast and slow times
is made explicit by a small parameter $\eps$,
and the quality of the approximation
can be measured in terms of~$\eps$.
By comparing the CSP manifold
with the slow manifold found in
Fenichel's
geometric singular perturbation theory
(GSPT, \cite{F-1979, J-1994}),
we showed that
each application of the CSP algorithm
improves the asymptotic accuracy of
the CSP manifold by one order of~$\eps$.

In this paper, 
we complete the analysis of
the CSP method
by focusing on the fast dynamics.
According to Fenichel's theory,
the fast--slow systems
we consider have,
besides a slow manifold,
a family of
fast stable fibers
along which
initial conditions 
tend toward
the slow manifold.
The base points of
these fibers 
lie on 
the slow manifold,
and
the dynamics 
near the slow manifold
can be decomposed into
a fast contracting component 
along the fast fibers
and
a slow component
governed by
the motion of
the base points
on the slow manifold.
By comparing
the CSP fibers
with the
tangent spaces of
the fast fibers at
their base points,
we show that
each application
of the CSP algorithm
also improves
the asymptotic accuracy
of the CSP fibers
by one order of~$\eps$.

Summarizing the results of~\cite{ZKK-2003}
and the present investigation,
we conclude that 
the CSP method provides
for the simultaneous approximation
of the slow manifold and 
the tangents to the fast fibers
at their base points.
If one is interested only
in the slow manifold,
then it suffices to implement
a reduced (one-step) version
of the algorithm.
On the other hand,
if one is interested
in both the slow and fast dynamics,
then it is necessary to use
the full (two-step) CSP algorithm.
Moreover, only the full CSP algorithm
allows for a linear matching
of any initial data with the dynamics
on the slow manifold.

This paper is organized as follows.
In Section~\ref{s-fast-slow},
we recall the relevant results
from Fenichel's theory
and set the framework for the CSP method.
In Section~\ref{s-CSP},
we outline the CSP algorithm
and state the main results:
Theorem~\ref{t-CSPMq}
concerning the approximation of the slow manifold,
which is a verbatim restatement
of~\cite[Theorem~3.1]{ZKK-2003};
and Theorem~\ref{t-CSPFq}
concerning the approximation of 
the tangent spaces of
the fast fibers.
The proof of Theorem~\ref{t-CSPFq}
is given in Section~\ref{s-Main}.
In Section~\ref{s-MMH model},
we revisit the 
Michaelis--Menten--Henri
mechanism of
enzyme kinetics
to illustrate
the CSP method 
and
the results of
this article.
Section~\ref{s-disc}
is devoted to a discussion of
methods for linearly projecting
initial conditions on
the slow manifold.

\section{Slow Manifolds and Fast Fibers \label{s-fast-slow}}
\setcounter{equation}{0}
Consider a general system of ODEs,
\be
  \frac{dx}{dt} = g(x) ,
  \Label{system}
\ee
for a vector-valued function
$x \equiv x(t) \in \mathbf{R}^{m+n}$
in a smooth vector field~$g$.
For the present analysis,
we assume that $n$ components of $x$
evolve on a time scale characterized
by the ``fast'' time $t$,
while the remaining $m$ components
evolve on a time scale characterized
by the ``slow'' time $\tau = \eps t$,
where $\eps$ is a small parameter.
(The explicit identification
of a small parameter $\eps$
is not necessary for
the applicability of the CSP method;
a separation of time scales is sufficient.)
We collect the slow variables in 
$y\in\mathbf{R}^m$
and the fast variables in
$z\in\mathbf{R}^n$.
Thus, the system~(\ref{system})
is equivalent to
either the ``fast system''
\begin{eqnarray}
  y' &=& \eps g_1 (y, z, \eps) ,  \Label{eq-y} \\
  z' &=& g_2 (y, z, \eps) , \Label{eq-z}
\end{eqnarray}
or the ``slow system''
\begin{eqnarray}
  \dot{y} &=& g_1 (y, z, \eps) ,  \Label{eq-y-slow} \\
  \eps \dot{z} &=& g_2 (y, z, \eps) . \Label{eq-z-slow}
\end{eqnarray}
(A prime~${}'$
denotes differentiation
with respect to~$t$,
a dot~$\dot{\ }$
differentiation
with respect to~$\tau$.)
The fast system is more appropriate
for the short-term dynamics,
the slow system for
the long-term dynamics
of the system~(\ref{system}).

In the limit as $\eps$ tends to 0,
the fast system
reduces formally to
a single equation
for the fast variable~$z$,
\be
    z' = g_2 (y,z,0) ,
\Label{eq-z-fastred}
\ee
where $y$ is a parameter,
while the slow system
reduces to
a differential equation
for the slow variable~$y$,
\be
    \dot{y} = g_1 (y,z,0) ,
\Label{eq-y-slowred}
\ee
with the algebraic constraint
$g_2 (y, z, 0) = 0$.

We assume that 
there exist a compact domain $K$
and a smooth function~$h_0$ defined on $K$
such that 
\be
    g_2 (y, h_0 (y), 0) = 0 , \quad y \in K .
\ee
The graph of $h_0$ defines
a critical manifold~$\Mo$,
\be
  \mathcal{M}_0
  =
  \{ (y, z) \in \mathbf{R}^{m+n} : z = h_0 (y) , \; y \in K \} ,
  \Label{M0}
\ee
and with each point $p = (y, h_0 (y)) \in \Mo$
is associated a fast fiber~$\mathcal{F}^p_0$,
\be
  \mathcal{F}^p_0
  =
  \{ (y, z) \in \mathbf{R}^{m+n} :
  z \in \mathbf{R}^n \} , \quad p \in \Mo .
  \Label{F0}
\ee
The points of $\Mo$ are fixed points
of Eq.~(\ref{eq-z-fastred}).
If the real parts of the
eigenvalues of $D_z g_2 (y, h_0(y), 0)$
are all negative, as we assume,
then $\Mo$ is asymptotically stable,
and
all solutions on~$\Fo^p$ 
contract exponentially toward~$p$.

If $\eps$ is positive but arbitrarily small,
Fenichel's theory~\cite{F-1979, J-1994}
guarantees that 
there exists a function $h_\eps$ 
whose graph is a slow manifold~$\Me$,
\be
  \Me = \{ (y,z)\in\mathbf{R}^{m+n} : z = h_\eps (y) , \; y \in K \} .
\Label{M-eps}
\ee
This manifold is
locally invariant
under the system dynamics,
and the dynamics on $\Me$
are governed by the equation
\be
  \dot{y} = g_1 ( y, h_\eps(y), \eps) ,
  \Label{traj}
\ee
as long as $y\in K$.
Fenichel's theory
also guarantees that
there exists
an invariant family~$\Fe$,
\be
    \Fe = \bigcup_{p\in\Me} \Fe^p ,
\Label{F-eps}
\ee
of fast stable fibers~$\Fe^p$
along which solutions relax to $\Me$.
The family is invariant
in the sense that,
if $\phi_t$ denotes the time-$t$ map
associated with Eq.~(\ref{system}),
then
\be
    \phi_t (\Fe^p) \subset \Fe^{\phi_t (p)} , \quad p \in \Me .
\ee
The collection of fast fibers $\Fe^p$
foliates a neighborhood of $\Me$.
Hence, the motion of
any point on $\Fe^p$ 
decomposes into
a fast contracting component
along the fiber
and a slow component 
governed by
the motion of the base point 
of the fiber.
Also,
$\Me$ is 
${\cal O}(\eps)$-close to~$\Mo$,
with
\be
  h_\eps (y)
  = h_0 (y) + \eps h_1 (y)
  + \eps^2 h_2 (y) + \cdots \,, \quad \eps \downarrow 0 ,
\Label{h-eps-exp}
\ee
and
$\Fe^p$ is
${\cal O}(\eps)$-close to $\Fo^p$
in any compact neighborhood of $\Me$.

\noindent\textbf{Remark 2.1.}
Typically, 
the manifold $\Me$ 
is not unique; 
there is a 
family of 
slow manifolds, 
all having
the same asymptotic expansion 
(\ref{h-eps-exp}) 
to all orders in $\eps$
but differing by 
exponentially small amounts
($\mathcal{O} (e^{-c/\eps}), c>0$ ).

\section{The CSP Method \label{s-CSP}}
\setcounter{equation}{0}
The CSP method focuses on
the dynamics of the vector field $g(x)$,
rather than on the dynamics
of the vector $x$ itself.

Writing a single differential
equation like~(\ref{system})
as a system of equations
amounts to choosing
a basis in the vector space.
For example,
in Eqs.~(\ref{eq-y})--(\ref{eq-z}),
the basis consists of the ordered set
of unit vectors in $\mathbf{R}^{m+n}$.
The coordinates of $g$
relative to this basis are
$\eps g_1$ and $g_2$.
If we collect the basis vectors
in a matrix in the usual way,
then we can express the relation between
$g$ and its coordinates in the form
\be
  g = \left(
  \begin{array}{cc}
  I_m & 0 \\
  0 & I_n
  \end{array} \right)
  \left(
  \begin{array}{c}
  \eps g_1 \\
  g_2
  \end{array} \right) .
\ee
Note that the basis chosen
for this representation
is the same at every point
of the phase space.
The CSP method is based on
a generalization of this idea,
where the basis is allowed to vary
from point to point,
so it can be tailored
to the local dynamics near $\Me$.

Suppose that we choose,
instead of a fixed basis,
a (point-dependent) basis~$A$
for~$\mathbf{R}^{m+n}$.
The relation between
the vector field $g$
and 
the vector $f$ of 
its coordinates 
relative to this
basis is
\be
  g = A f .
\Label{g=Af}
\ee
Conversely,
\be
  f = B g ,
\Label{f=Bg}
\ee
where $B$ is
the left inverse of $A$,
$BA = I$ on 
$\mathbf{R}^{m+n}$. 
In the convention of the CSP method,
$A$ is a matrix of column vectors
(vectors in $\mathbf{R}^{m+n}$)
and $B$ a matrix of row vectors
(functionals on $\mathbf{R}^{m+n}$).

The CSP method focuses
on the dynamics of the vector $f$.
Along a trajectory of the system~(\ref{system}),
$f$ satisfies the ODE
\be
  \frac{df}{dt}
  = 
  \Lambda f ,
\Label{f-ODE}
\ee
where $\Lambda$
is a linear
operator~\cite{M-1996, ZKK-2003},
\be
  \Lambda
  =
  B (Dg) A + \frac{dB}{dt} A
  =
  B (Dg) A - B \frac{dA}{dt}
  =
  B [A, g] .
  \Label{Lambda}
\ee
Here, $Dg$ is the Jacobian of $g$,
$dB/dt = (DB) g$, $dA/dt = (DA) g$,
and $[A,g]$ is the Lie bracket of $A$
(taken column by column) and $g$.
The Lie bracket of any two vectors
$a$ and $g$ is
$[a, g] = (Dg)a - (Da)g$;
see~\cite{O-1986}.

It is clear from Eq.~(\ref{f-ODE})
that the dynamics of $f$
are governed by $\Lambda$,
so the CSP method focuses
on the structure of $\Lambda$.

\noindent\textbf{Remark 3.1.}
It is useful to see
how $\Lambda$ transforms
under a change of basis.
If $C$ is an invertible square matrix
representing a coordinate transformation
in $\mathbf{R}^{m+n}$, and
$\hat{A} = AC$ and
$\hat{B} = C^{-1} B$,
then
\begin{eqnarray}
  \hat{\Lambda}
  &=&
  \hat{B} (Dg) \hat{A} - \hat{B} \frac{d\hat{A}}{dt}
  =
  C^{-1} B (Dg) AC - C^{-1} B \frac{d(AC)}{dt} \nonumber \\
  &=&
  C^{-1} B (Dg) AC
  - C^{-1} B \left( \frac{dA}{dt} C + A \frac{dC}{dt} \right) \nonumber \\
  &=&
  C^{-1} \Lambda C - C^{-1} \frac{dC}{dt} .
\Label{L-transf}
\end{eqnarray}
Hence, $\Lambda$ does not
transform as a matrix,
unless $C$ is constant.

\subsection{Decompositions \label{ss-decomp}}
Our goal is to decompose
the vector $f$ into
its fast and slow components.
Suppose, therefore, that
we have a decomposition of
this type,
$f =
\left( \begin{array}{c} 
      f^1 \\ f^2 
       \end{array}
\right)$,
where 
$f^1$ and $f^2$ 
are of length
$n$ and $m$, respectively,
but not necessarily
fast and slow everywhere. 
The decomposition suggests
corresponding decompositions
of the matrices $A$ and $B$,
namely
$A =
  \left( A_1, A_2 \right)$
and
$B =
\left( \begin{array}{c} 
       B^1 \\ B^2 
       \end{array}
\right)$,
where
$A_1$ is an $(m+n) \times n$ matrix,
$A_2$ an $(m+n) \times m$ matrix,
$B^1$ an $n \times (m+n)$ matrix,
and
$B^2$ an $m \times (m+n)$ matrix.
Then,
$f^1 = B^1 g$ and $f^2 = B^2 g$.

The decompositions of $A$ and $B$ lead,
in turn, to a decomposition of $\Lambda$,
\be
  \Lambda
  =
  \left(
  \begin{array} {cc}
  \Lambda^{11} & \Lambda^{12} \\
  \Lambda^{21} & \Lambda^{22}
  \end{array} \right)
  =
  \left(
  \begin{array} {cc}
  B^1 [A_1, g] & B^1 [A_2, g] \\
  B^2 [A_1, g] & B^2 [A_2, g]
  \end{array} \right) .
\Label{L-blocks}
\ee
The off-diagonal blocks
$\Lambda^{12}$
and
$\Lambda^{21}$
are, in general, not zero,
so the equations governing the
evolution of the coordinates
$f^1$ and $f^2$ are coupled.
Consequently,
$f^1$ and $f^2$ cannot be identified
with the fast and slow coordinates of $g$
globally along trajectories.
The objective of the CSP method is
to construct local coordinate systems
(that is, matrices $A$ and $B$)
that lead to a block-diagonal structure
of $\Lambda$.
We will see,
in the next section,
that such a structure is
associated with
a decomposition
in terms of
the slow manifold 
and
the fast fibers.

\noindent\textbf{Remark 3.2.}
Note that the identity $BA=I$
on $\mathbf{R}^{m+n}$
implies four identities,
which are summarized in
the matrix identity
\be
 \left(
 \begin{array}{cc}
  B^1 A_1 & B^1 A_2 \\
  B^2 A_1 & B^2 A_2
  \end{array} \right)
 =
 \left(
 \begin{array}{cc}
 I_n & 0 \\
 0 & I_m
 \end{array} \right) .
\Label{BA}
\ee

\subsection{Block-Diagonalization of $\Lambda$ \label{ss-blockL}}
In this section
we analyze the properties of $\Lambda$ 
relative to a fast--slow decomposition
of the dynamics near $\Me$.

Let
${\mathcal T}_p\Fe$
and
${\mathcal T}_p\Me$
denote the tangent spaces
to the fast fiber
and the slow manifold,
respectively,
at the base point $p$
of the fiber on $\Me$.
(Note that
dim${\mathcal T}_p\Fe=n$
and
dim${\mathcal T}_p\Me=m$.)
These two linear spaces
intersect transversally,
because $\Me$ is normally hyperbolic
and compact, so
\be
    \mathbf{R}^{m+n} 
= 
    \mathcal{T}_p \Fe \oplus \mathcal{T}_p\Me ,
\quad p\in\Me .
\Label{Tsplit}
\ee
Let $A_f$ be 
an $(m+n) \times n$
matrix 
whose columns form 
a basis
for ${\mathcal T}_p\Fe$
and $A_s$
an $(m+n) \times m$
matrix 
whose columns form 
a basis
for ${\mathcal T}_p\Me$,
and let
$A = (A_f, A_s)$.
(We omit the subscript $p$.)
Then $A$ is a (point-dependent) basis
for $\mathbf{R}^{m+n}$
that respects the
decomposition~(\ref{Tsplit}).
We recall that
${\cal T}\Me \equiv \bigcup_{p\in\Me} (p,{\cal T}_p\Me)$
and
${\cal T}\Fe \equiv \bigcup_{p\in\Me} (p,{\cal T}_p\Fe)$
are the tangent bundles of
the slow manifold 
and 
the family of the fast fibers,
respectively.
(A general treatment of
tangent bundles of manifolds
is given in
\cite[Section~1.7]{DFN-1985}.)

The decomposition~(\ref{Tsplit})
induces a dual decomposition,
\be
    \mathbf{R}^{m+n} = \mathcal{N}_p \Me \oplus \mathcal{N}_ p\Fe ,
\quad p\in\Me ,
\Label{Nsplit}
\ee
where
${\mathcal N}_p\Me$
and
${\mathcal N}_p\Fe$
are the duals
of
${\mathcal T}_p\Me$ 
and
${\mathcal T}_p\Fe$,
respectively,
in $\mathbf{R}^{m+n}$.
(Note that
dim${\mathcal N}_p\Me=n$
and
dim${\mathcal N}_p\Fe=m$.)
The corresponding decomposition
of $B$ is
$B
=
\left(\begin{array}{c}
B^{s \perp}\\ 
B^{f \perp}
\end{array}\right)$,
where the rows of $B^{s \perp}$
form a basis for $\mathcal{N}_p \Me$
and the rows of $B^{f \perp}$
a basis for $\mathcal{N}_p \Fe$.
Furthermore,
\be
    \left(
          \begin{array}{cc}
          B^{s \perp} A_f & B^{s \perp} A_s \\
          B^{f \perp} A_f & B^{f \perp} A_s \\
          \end{array}
    \right)
    =
    \left(
           \begin{array}{cc}
           I_n & 0\\
           0 & I_m
           \end{array}
    \right) .
\ee

The decompositions of $A$ and $B$ lead,
in turn, to a decomposition of $\Lambda$,
\be
  \Lambda
  =
  \left(
  \begin{array} {cc}
  B^{s \perp} [A_f, g] & B^{s \perp} [A_s, g] \\
  B^{f \perp} [A_f, g] & B^{f \perp} [A_s, g]
  \end{array} \right) .
\Label{L-fsblocks}
\ee
This decomposition is similar to,
but different from,
the decomposition~(\ref{L-blocks}).
The following lemma 
shows that 
its off-diagonal blocks 
are zero.

\begin{LEMMA} \label{l-Ldiag}
The off-diagonal blocks
in the representation~(\ref{L-fsblocks})
of $\Lambda$
are zero at each point $p \in \Me$.
\end{LEMMA}

\begin{PROOF}
Since
$B^{s \perp} A_s = 0 $
on $\Me$
and
$\Me$ is invariant,
we have
\be
    \frac{d}{dt} \left(B^{s \perp} A_s\right)
=
    D(B^{s \perp} A_s) g
=
      (DB^{s \perp}) (g,A_s)
    + B^{s \perp} ((DA_s) g)
=
    0 .
\Label{DBAg=0}
\ee
($DB^{s \perp}$ is a symmetric bilinear form;
its action on a matrix must be understood
as column-wise action.)

Also,
$g \in {\mathcal T} \Me$,
so
$B^{s \perp} g = 0$ 
on $\Me$.
Hence,
the directional derivative 
along $A_s$
(taken column by column)
at points on
$\Me$
also vanishes,
\be
    D(B^{s \perp} g) A_s
=
      (DB^{s \perp}) (A_s,g)
    + B^{s \perp} (Dg) A_s
=
    0 .
\Label{DBgA=0}
\ee
Subtracting 
Eq.~(\ref{DBAg=0})
from
Eq.~(\ref{DBgA=0}),
we obtain the identity
\be
    B^{s \perp} [A_s,g] 
=
    B^{s \perp} \left(\left(Dg\right) A_s - (DA_s) g\right) 
=  
    0 .
\Label{B[g,A]=0}
\ee

The proof for
the lower left block
is more involved,
since
the fast fibers
are invariant as
a family.
Assume that
the fiber $\Fe^p$
at $p \in \Me$
is given implicitly
by the equation
$F(q;p) = 0$, $q \in \Fe^p$.
Then
the rows of
$(D_q F) (q;p)$
form
a basis for
${\cal N}_q \Fe$,
so there exists an
invertible matrix $C$ such that
$B^{f \perp} = C (D_q F)$.

Since 
the rows of
$(D_q F) (q;p)$ span
${\cal N}_q \Fe$,
we have
$(D_q F) (q;p) A_f (q) = 0$.
This identity holds,
in particular,
along solutions
of~(\ref{system}),
so
\begin{eqnarray}
    \frac{d}{dt} ((D_q F) (q;p) A_f (q))
&=&
      \left((D_q^2 F) (q;p) \right) (g(q), A_f(q)) \nonumber \\
&&\mbox{}
    + \left((D_{pq}F) (q;p) \right) (g(p), A_f(q)) \nonumber\\
&&\mbox{}
    + \left((D_qF) (q;p) \right) \left(DA_f(q)\right) g(q) \nonumber\\
&=&
    0 .
\Label{DFAg=0}
\end{eqnarray}
The family of
the fast fibers
is invariant
under
the flow
associated with
(\ref{system}),
so
if $F(q;p) = 0$,
then also
$F(q(t);p(t)) = 0$
and, hence,
\begin{eqnarray}
    \frac{d F(q;p)}{dt}
=
      \left((D_qF) (q;p) \right) g(q)
    + \left((D_pF) (q;p) \right) g(p)
=
    0 .
\Label{DFgA=0-aux}
\end{eqnarray}
Next, 
we take 
the directional derivative of 
both members of 
this equation
along $A_f$,
keeping in mind that
$(Dg)(p) A_f(q) = 0$
because
the base point $p$
does not vary
along $A_f$.
(Recall that
the columns of
$A_f(q)$
span
${\cal T}_q\Fe$.)
We find
\begin{eqnarray}
&&\mbox{}
      \left((D_q^2F) (q;p) \right) \left(A_f(q),g(q)\right)
    +  \left((D_qF) (q;p) \right) \left(Dg(q)\right) A_f(q)
\nonumber\\
&&\mbox{}
    + \left((D_{pq}F) (q;p) \right) (A_f(q),g(p))
=
    0 .
\Label{DFgA=0}
\end{eqnarray}
But
the bilinear forms
$D_q^2 F$
and
$D_{pq}F$
are symmetric,
so subtracting
Eq.~(\ref{DFAg=0})
from Eq.~(\ref{DFgA=0})
and letting $q = p$,
we obtain the identity
\begin{eqnarray}
  (D_qF) (p;p) \left(\left(Dg\right) A_f - (DA_f) g\right)(p) 
=
  0 .
\Label{F[A,g]=0}
\end{eqnarray}
Hence,
$B^{f \perp} [A_f,g](p)
=
C (D_qF) (p;p) [A_f,g](p)
=
0$,
and 
the proof of the lemma is complete.
\end{PROOF}

The lemma implies that
the representation~(\ref{L-fsblocks})
is  block-diagonal,
\be
  \Lambda
  =
  \left(
  \begin{array} {cc}
  B^{s \perp} [A_f, g] & 0 \\
  0 & B^{f \perp} [A_s, g]
  \end{array} \right) .
\Label{L-diag}
\ee
Consequently,
the decomposition~(\ref{Tsplit})
reduces $\Lambda$.
In summary, 
if we can construct bases
$A_f$ and $A_s$,
then we will have achieved
a representation of $\Lambda$
where the fast and slow components
remain separated at all times
and the designation of
fast and slow
takes on a global meaning.

\subsection{The CSP Algorithm \label{ss-CSP}}
The CSP method is
a constructive algorithm
to approximate $A_f$ and $A_s$.
One typically initializes
the algorithm with
a constant matrix 
$\Af^{(0)}$,
\be
  \Af^{(0)}
  =
  \left( \Af_1^{(0)}, \, \Af_2^{(0)} \right)
  =
  \left( 
  \begin{array}{cc}
    \Af_{11}^{(0)} & \Af_{12}^{(0)} \\
    \Af_{21}^{(0)} & \Af_{22}^{(0)}
  \end{array}
  \right) .
  \Label{A(0)-full}
\ee
Here,
$\Af_{11}^{(0)}$ is an $m \times n$ matrix,
$\Af_{22}^{(0)}$ an $n \times m$ matrix,
and the off-diagonal blocks
$\Af_{12}^{(0)}$ and $\Af_{21}^{(0)}$
are full-rank square matrices of order
$m$ and $n$, respectively.
A common choice is
$\Af_{11}^{(0)} = 0$.
We follow this convention and assume,
henceforth, that $\Af_{11}^{(0)} = 0$,
\be
  \Af^{(0)}
  =
  \left( \Af_1^{(0)}, \, \Af_2^{(0)} \right)
  =
  \left( 
  \begin{array}{cc}
    0            & \Af_{12}^{(0)} \\
    \Af_{21}^{(0)} & \Af_{22}^{(0)}
  \end{array}
  \right) .
  \Label{A(0)}
\ee
(Other choices are discussed in \cite{ZKK-2003}.)
The left inverse of $\Af^{(0)}$ is
\begin{eqnarray}
  \Bf_{(0)}
&=& 
  \left(
  \begin{array}{c} \Bf^1_{(0)} \\ \Bf^2_{(0)} \end{array}
  \right)
  =
  \left( 
  \begin{array}{cc}
    \Bf^{11}_{(0)} & \Bf^{12}_{(0)} \\
    \Bf^{21}_{(0)} & 0
  \end{array}
  \right) \nonumber \\
  &=& 
  \left(
  \begin{array}{cc}
    - ( \Af_{21}^{(0)} )^{-1}
    \Af_{22}^{(0)}
    ( \Af_{12}^{(0)} )^{-1}
    & ( \Af_{21}^{(0)} )^{-1} \\
    ( \Af_{12}^{(0)} )^{-1}
    & 0
  \end{array}
  \right) .
  \Label{B(0)}
\end{eqnarray}

The algorithm proceeds iteratively.
For $q = 0, 1, \ldots\,$,
one first defines the operator
$\bfL_{(q)}$
in accordance with
Eq.~(\ref{Lambda}),
\be
  \bfL_{(q)}
  =
  \Bf_{(q)} (Dg) \Af^{(q)} - \Bf_{(q)} \frac{d\Af^{(q)}}{dt}
  =
  \left(
  \begin{array}{cc}
    \bfL_{(q)}^{11} & \bfL_{(q)}^{12} \\
    \bfL_{(q)}^{21} & \bfL_{(q)}^{22}
  \end{array}
  \right) ,
  \Label{Lambda(q)}
\ee
and matrices
$\Uf_{(q)}$
and
$\Lf_{(q)}$,
\be
  \Uf_{(q)}
  =
  \left(
  \begin{array}{cc}
    0 & (\bfL_{(q)}^{11})^{-1} \bfL_{(q)}^{12} \\
    0 & 0
  \end{array}
  \right) , \quad
  \Lf_{(q)}
  =
  \left(
  \begin{array}{cc}
    0 & 0 \\
    \bfL_{(q)}^{21} (\bfL_{(q)}^{11})^{-1} & 0
  \end{array}
  \right) .
\Label{U(q),L(q)}
\ee
Then one updates
$A^{(q)}$ and $B_{(q)}$
according to the formulas
\begin{eqnarray}
  \Af^{(q+1)}
  &=&
  \Af^{(q)} ( I - \Uf_{(q)} ) ( I + \Lf_{(q)} ) , \Label{A(q+1)} \\
  \Bf_{(q+1)}
  &=&
  ( I - \Lf_{(q)} ) ( I + \Uf_{(q)} ) \Bf_{(q)} ,  \Label{B(q+1)}
\end{eqnarray}
and returns to Eq.~(\ref{Lambda(q)})
for the next iteration.

\noindent\textbf{Remark 3.3.}
Lam and Goussis~\cite{L-1993}
perform the update~(\ref{A(q+1)})--(\ref{B(q+1)}) 
in two steps.
The first step 
corresponds to 
the postmultiplication of 
$\Af^{(q)}$ with $I - \Uf_{(q)}$ 
and
premultiplication of 
$\Bf_{(q)}$ with $I + \Uf_{(q)}$,
the second step to 
the subsequent postmultiplication of 
$\Af^{(q)}(I - \Uf_{(q)})$ with $I + \Lf_{(q)}$ 
and
premultiplication of 
$(I + \Uf_{(q)})\Bf_{(q)}$ with $I - \Lf_{(q)}$.

\subsection{Approximation of the Slow Manifold \label{ss-slow}}
After $q$ iterations,
the CSP condition
\be
  B_{(q)}^1 g = 0 , \quad q = 0, 1, \ldots \,,
\Label{CSPM-q}
\ee
identifies those points
where the fast amplitudes
vanish with respect to
the then current basis.
These points define a manifold
that is an approximation
for the slow manifold $\Me$.

For $q=0$, $B_{(0)}^1$ is constant
and given by Eq.~(\ref{B(0)}).
Hence, the CSP condition~(\ref{CSPM-q})
reduces to the constraint
$g_2 (y, z, \eps) = 0$.
In general,
this constraint is satisfied by
a function $z = \psi_{(0)}(y,\eps)$.
The graph of this function
defines~$\Kf^{(0)}$,
the CSP manifold (CSPM) of order zero.
Since the constraint reduces
at leading order to the equation
$g_2 (y, z, 0) = 0$,
which is satisfied by 
the function $z = h_0 (y)$,
$\Kf^{(0)}$
may be chosen to coincide 
with $\Mo$ to leading order;
see Eq.~(\ref{M0}).

For $q = 1, 2, \ldots\,$,
the CSP condition takes the form
\be
  B^1_{(q)}(y,\psi_{(q-1)}(y,\eps),\eps)
  g(y,z,\eps) = 0 , \quad
  q = 1, 2, \ldots \,.
\ee
The condition
is satisfied by a function
$z = \psi_{(q)}(y,\eps)$,
and the manifold
\be
  \Kf^{(q)}
  = \{ (y,z): z = \psi_{(q)} (y, \eps) , \; y \in K \} ,
  \quad q = 0, 1, \ldots
\Label{L-q}
\ee
defines the CSP manifold (CSPM) of order $q$,
which is an approximation of $\Me$.
The following theorem regarding
the quality of the approximation
was proven in~\cite{ZKK-2003}.

\begin{THEOREM}
\label{t-CSPMq}
\cite[Theorem~3.1]{ZKK-2003}
The asymptotic expansions of
the CSP manifold 
$\Kf^{(q)}$
and 
the slow manifold
$\Me$
agree
up to and including terms
of $\mathcal{O} (\eps^q)$,
\be
\psi_{(q)} ( \cdot\,, \eps)
=
\sum_{j=0}^q \eps^j h_j + \mathcal{O} (\eps^{q+1}) , \quad
\eps \downarrow 0 , \quad
  q = 0, 1, \ldots\,.
\ee
\end{THEOREM}

\subsection{Approximation of the Fast Fibers \label{ss-fast}}
We now turn our attention
to the fast fibers.
The columns
of $A_f(y,h_\eps(y))$ span
the tangent space
to the fast fiber
with base point $p = (y,h_\eps(y))$,
so we expect that
$A_1^{(q)}$ defines
an approximation
for the same space 
after $q$ applications
of the CSP algorithm.
We denote this approximation
by $\mathcal{L}_\eps^{(q)}(y)$
and refer to it as the
CSP fiber (CSPF) of order $q$
at $p$,
\be
  \mathcal{L}_\eps^{(q)} (y)
  =
  \mbox{span (cols } (A_1^{(q)} (y,\psi_{(q)}(y,\eps),\eps))) .
\Label{CSPF-def}
\ee
We will shortly estimate
the asymptotic accuracy
of the approximation,
but before doing so
we need to make an
important observation.

Each application of
the CSP algorithm
involves two steps,
see Remark~3.3.
The first step involves $U$
and serves to push
the order of magnitude
of the upper right block
of $\Lambda$
up by one,
the second step
involves $L$ and serves
the same purpose
for the lower left block.
The two steps are consecutive.
At the first step of the $q$th iteration,
one evaluates
$B^1_{(q)}$
on
$\mathcal{K}_\eps^{(q-1)}$
to find
${\cal K}_\eps^{(q)}$
by solving the CSP
condition~(\ref{CSPM-q})
for the function
$\psi_{(q)}$.
One then uses
this expression
in the second step
to update $A$ and $B$,
thus effectively
evaluating $A_1^{(q)}$
on ${\cal K}_\eps^{(q)}$
rather than on
${\cal K}_\eps^{(q-1)}$.

The following theorem 
contains our main result.

\begin{THEOREM}
\label{t-CSPFq}
The asymptotic expansions of
${\cal L}_\eps^{(q)}(y)$
and
${\cal T}_p \Fe$,
where $p = (y,h_\eps(y)) \in \Me$,
agree up to and including terms of
${\cal O}(\eps^q)$,
for all $y\in K$
and
for $q = 0, 1, \ldots\,$.
\end{THEOREM}

Theorem~\ref{t-CSPFq}
implies that
the family
${\cal L}_\eps^{(q)} \equiv \bigcup_{p\in\Me} (p,{\cal L}_\eps^{(q)}(y))$
is an 
${\cal O}(\eps^q)$-approximation to
the tangent bundle
${\cal T}\Fe$.

The proof of 
Theorem~\ref{t-CSPFq}
is given in 
Section~\ref{s-Main}.
The essential idea 
is to show that,
at each iteration,
the asymptotic order
of the off-diagonal blocks
of $\Lambda_{(q)}$
increases by one
and
$A_1^{(q)}$
and
$B^2_{(q)}$
become
fast
and
fast$^\perp$,
respectively,
to one higher order.
As a consequence,
in the limit as
$q \to \infty$,
$\Lambda_{(q)} \to \Lambda$,
$A^{(q)} \to A$,
and
$B_{(q)} \to B$,
where 
$\Lambda$, $A$, and $B$
are ideal
in the sense 
described in 
Section~\ref{ss-blockL}.

\noindent\textbf{Remark 3.4.}
If,
in the second step of
the CSP algorithm,
$A_1^{(q)}$ were evaluated on
$\Kf^{(q-1)}$
instead of on
$\Kf^{(q)}$,
the approximation of ${\cal T} \Fe$
might be only
${\cal O}(\eps^{q-1})$-accurate.
However,
see Section~\ref{s-MMH model}
for an example where 
the approximation is still
${\cal O}(\eps^{q})$.


\section{Proof of Theorem~\ref{t-CSPFq} \label{s-Main}}
\setcounter{equation}{0}
The proof of Theorem~\ref{t-CSPFq}
is by induction on $q$.
Section~\ref{ss-L}
contains an auxiliary
lemma that shows that
each successive application
of the CSP algorithm
pushes $\Lambda$
closer to
block-diagonal form.
The induction hypothesis
is formulated in
Section~\ref{ss-ind hypo},
the hypothesis is shown
to be true for $q=0$
in Section~\ref{ss-q=0},
and the induction step
is taken in
Section~\ref{ss-q=1,2,...}.

\subsection{Asymptotic Estimates of $\bfL$ \label{ss-L}}
As stated in
Section~\ref{s-CSP},
the goal of
the CSP method is to
reduce $\bfL$ to
block-diagonal form.
This goal is approached
by the repeated application
of a two-step algorithm.
As shown in~\cite{ZKK-2003},
the first step of the algorithm
is engineered
so that each application
increases the
asymptotic accuracy
of the upper-right block 
$\bfL_{(q)}^{12}$
by one order of $\eps$;
in particular,
$\bfL^{12}_{(q)} = {\cal O}(\eps^q)$
on $\Kf^{(q)}$
\cite[Eq.~(5.25)]{ZKK-2003}. 
We now complete
the picture
and show that
each application of
the second step
increases the
asymptotic accuracy
of the lower-left block
$\bfL^{21}_{(q)}$
by one order of $\eps$,
when the information obtained in
the first step of
the same iteration
is used.
In particular,
$\bfL^{21}_{(q)} = {\cal O}(\eps^{q+1})$
on $\Kf^{(q+1)}$,
where
$\Kf^{(q+1)}$ has been obtained in
the first step of
the $(q+1)$th refinement.

\begin{LEMMA}
\label{l-L-estim}
For
$q = 0,1,\ldots$,
\be
    \bfL_{(q)} = \left(
                 \begin{array}{cc}
                 \bfL^{11}_{(0,0)} + {\cal O}(\eps)
                 & \eps^q \bfL^{12}_{(q,q)}\\
                   \eps^{q+1} \bfL^{21}_{(q,q+1)}
                 & \eps \bfL^{22}_{(1,1)}+{\cal O}(\eps^2)
                 \end{array}
               \right) ,
\ee
when
$\bfL_{(q)}$ is  
evaluated on $\Kf^{(q+1)}$.
\end{LEMMA}
\begin{PROOF}
The proof is
by induction.
The desired 
estimates of
$\bfL^{11}_{(q)}$, 
$\bfL^{12}_{(q)}$, 
and
$\bfL^{22}_{(q)}$
on
$\Kf^{(q)}$
were established
in~\cite[Eqs.~(5.24), (5.25), (5.27)]{ZKK-2003}.
Since the asymptotic expansions of
$\Kf^{(q+1)}$
and
$\Kf^{(q)}$
differ only at
terms of
${\cal O}(\eps^{q+1})$
or higher
(\cite[Theorem~3.1]{ZKK-2003}),
these estimates of
$\bfL^{11}_{(q)}$, 
$\bfL^{12}_{(q)}$, 
and
$\bfL^{22}_{(q)}$
are true also on
$\Kf^{(q+1)}$.
It only remains to 
estimate $\bfL^{21}_{(q)}$.

Consider the case $q=0$.
Let $\bfL^{21}_{(0,j)}$
be
the coefficient of
$\eps^j$
in the asymptotic expansion of
$\bfL^{21}_{(0)} (y,\psi_{(1)}(y),\eps)$.
The estimate
$\bfL^{21}_{(0)} = {\cal O}(\eps)$
on $\Kf^{(1)}$
follows if we can show that
$\bfL^{21}_{(0,0)} = 0$.
It is already stated in
\cite[Eq.~(4.30)]{ZKK-2003}
that
$\bfL^{21}_{(0,0)} = 0$
on $\Kf^{(0)}$.
Furthermore,
\cite[Theorem~3.1]{ZKK-2003}
implies that
the asymptotic expansions of
$\psi_{(1)}$
and
$\psi_{(0)}$
agree to
leading order.
Thus,
the asymptotic expansions of
$\bfL^{21}_{(0)} (y,\psi_{(0)}(y),\eps)$
and
$\bfL^{21}_{(0)} (y,\psi_{(1)}(y),\eps)$
also agree
to leading order,
and
the result follows.

Now, assume that
the asymptotic estimate holds
for $0,1,\ldots\,,q$.
From Eq.~(\ref{L-transf})
we obtain
\begin{eqnarray}
    \bfL^{21}_{(q+1)} 
&=&
     \bfL^{21}_{(q)}
    - L_{(q)} \bfL^{11}_{(q)}   
    + \bfL^{22}_{(q)} L_{(q)}
    - L_{(q)} \bfL^{12}_{(q)} L_{(q)}
    - \bfL^{21}_{(q)} U_{(q)} L_{(q)}
\nonumber\\
&&\mbox{}
    - L_{(q)} U_{(q)} \bfL^{21}_{(q)}
    + L_{(q)} \bfL^{11}_{(q)} U_{(q)} L_{(q)}
    - L_{(q)} U_{(q)} \bfL^{22}_{(q)} L_{(q)}
\nonumber\\
&&\mbox{}
    + L_{(q)} U_{(q)} \bfL^{21}_{(q)} U_{(q)} L_{(q)}
    + \left(DL_{(q)}\right) g
    + L_{(q)} \left(\left(DU_{(q)}\right) g\right) L_{(q)} .~~~
\Label{L21-transf}
\end{eqnarray} 
The first two terms in 
the right member
sum to zero,
by virtue of
the definition~(\ref{U(q),L(q)})
of $L_{(q)}$.
The next seven terms
are all
${\cal O}(\eps^{q+2})$
or higher,
by virtue of
the induction hypothesis.
Finally,
the last two terms
are also
${\cal O}(\eps^{q+2})$ 
or higher,
by 
the induction hypothesis
and
\cite[Lemma~A.2]{ZKK-2003}.
\end{PROOF}

\subsection{The Induction Hypothesis \label{ss-ind hypo}}
The CSPF of order $q$,
${\mathcal L}_\eps^{(q)}(y)$,
is defined in 
Eq.~(\ref{CSPF-def})
to be
the linear space
spanned by
the columns of
the fast component,
$A^{(q)}_1(y,\psi_{(q)},\eps)$, of
the basis $A^{(q)}$.
Thus, 
to prove Theorem~\ref{t-CSPFq},
it suffices to 
show that
the asymptotic expansions of
$A^{(q)}_1(y,\psi_{(q)},\eps)$
and
the space 
tangent to
the fast fiber,
${\cal T}_p\Fe$,
agree up to and including terms of
${\cal O}(\eps^{q})$,
for
$p = (y,h_\eps(y))$
and
for $q = 0, 1, \ldots\,$.
The central idea of
the proof
is to 
show that
each successive
application of the CSP method
pushes
the projection of 
$A_1^{(q)}$ on
${\cal T}\Me$
along ${\cal T}\Fe$
to one higher order in $\eps$.

We express $A^{(q)}$,
generated after
$q$ applications of
the CSP algorithm,
in terms of
the basis $A$,
\be
    A^{(q)} (y,z,\eps) 
= 
    A (y,h_\eps,\eps) 
    Q^{(q)} (y,z,\eps) , \quad
  q = 0,1,\ldots \,.
\Label{A-split}
\ee
Since
$B_{(q)}$ and $B$
are the left inverses of
$A^{(q)}$ and $A$,
respectively,
we also have
\be
    B_{(q)} (y,z,\eps) 
= 
    R_{(q)} (y,z,\eps) 
    B (y,h_\eps,\eps) , \quad
  q = 0,1,\ldots\, ,
\Label{B-split}
\ee
where 
$R_{(q)} \equiv (Q^{(q)})^{-1}$.
Introducing the block structure of
$Q^{(q)}$ and $R_{(q)}$,
\begin{eqnarray}
    Q^{(q)} = \left(
                 \begin{array}{cc}
                 Q^{(q)}_{1f} & Q^{(q)}_{2f}\\
                 Q^{(q)}_{1s} & Q^{(q)}_{2s}
                 \end{array}
             \right) ,
\quad
    R_{(q)} = \left(
                 \begin{array}{cc}
                 R_{(q)}^{1s\perp} & R_{(q)}^{1f\perp}\\
                 R_{(q)}^{2s\perp} & R_{(q)}^{2f\perp}
                 \end{array}
             \right) ,
\Label{QR-def}
\end{eqnarray}
we rewrite Eqs.~(\ref{A-split}) 
and (\ref{B-split}) as
\begin{eqnarray}
    A^{(q)}_1 = A_f Q^{(q)}_{1f} 
                + A_s Q^{(q)}_{1s} ,  \quad
    A^{(q)}_2 = A_f Q^{(q)}_{2f} 
                + A_s Q^{(q)}_{2s} ,
\Label{Ai-split}
\end{eqnarray}
and
\begin{eqnarray}
    B^1_{(q)} = R_{(q)}^{1s\perp} B^{s\perp} 
                + R_{(q)}^{1f\perp} B^{f\perp} , \quad
    B^2_{(q)} = R_{(q)}^{2s\perp} B^{s\perp} 
                + R_{(q)}^{2f\perp} B^{f\perp} , 
\Label{Bi-split}
\end{eqnarray}
for $q = 0,1,\ldots\,$.

Equation~(\ref{Bi-split})
shows that
$A_s Q^{(q)}_{1s} $ is
the projection of
$A_1^{(q)}$ on 
${\cal T}\Me$.
Thus,
to establish Theorem~\ref{t-CSPFq},
we only need to
prove the asymptotic estimate
$Q^{(q)}_{1s} = {\cal O}(\eps^{q+1})$.
The proof is
by induction on
$q$,
where
the induction hypothesis
is
\begin{eqnarray}
    Q^{(q)} (\cdot\,,\psi_{(q)},\eps)
&=&
    \left( 
        \begin{array}{cc}
        {\cal O}(1)
        & {\cal O}(\eps^q)\\
        {\cal O}(\eps^{q+1})
        & {\cal O}(1)
        \end{array}
    \right) ,
\Label{Q(q)}\\
    R_{(q)} (\cdot\,,\psi_{(q)},\eps)
&=&
    \left( 
        \begin{array}{cc}
        {\cal O}(1)
        & {\cal O}(\eps^q)\\
        {\cal O}(\eps^{q+1})
        & {\cal O}(1)
        \end{array}
    \right) ,
\quad q=0,1,\ldots\, .
\Label{R(q)}
\end{eqnarray}

\noindent\textbf{Remark 4.1.}
Although
the estimate of
$Q^{(q)}_{1s}$
is sufficient to
establish Theorem~\ref{t-CSPFq},
we provide
the estimates of
all the blocks
in Eqs.~(\ref{Q(q)})--(\ref{R(q)})
because they will 
be required in
the induction step.

The validity of
Eqs.~(\ref{Q(q)})--(\ref{R(q)})
for $q = 0$
is shown in
Section~\ref{ss-q=0}.
The induction step is
carried out in
Section~\ref{ss-q=1,2,...}.

\subsection{Proof of Theorem~\ref{t-CSPFq} for $q=0$ \label{ss-q=0}}
We fix
$q=0$
and verify
the induction hypothesis
for $Q^{(0)}$ and $R_{(0)}$.
By Eq.~(\ref{A-split})
\be
    Q^{(0)} = B A^{(0)} ,
\Label{Q0}
\ee
whence
\begin{eqnarray}
    Q^{(0)}
&=& 
    \left(
           \begin{array}{cc}
               B^{s\perp} A^{(0)}_1 
             & B^{s\perp} A^{(0)}_2 \\
               B^{f\perp} A^{(0)}_1 
             & B^{f\perp} A^{(0)}_2 
           \end{array} 
    \right) .
\Label{blox-Q0}
\end{eqnarray}
It suffices to 
show that
the lower-left block
is zero
to leading order,
since the other blocks
are all ${\cal O}(1)$.
We do this
by showing that
$Q^{(0,0)}_{1s} = 0$.
By Eq.~(\ref{blox-Q0}),
\be
    Q^{(0,0)}_{1s}
=
    B^{f\perp}_0 A^{(0,0)}_1 .
\Label{Q1s(0,0)}
\ee
$B^{f\perp}_0$ spans
${\cal N}_p\Fo$
for every $p\in{\mathcal K}_\eps^{(0)}$.
Also,
$z$ is constant
on ${\cal N}_p\Fo$,
so
$B^{f\perp}_0
=
(B^{1f\perp},0)$,
where $B^{1f\perp}$ is
a full-rank matrix of 
size $m$.
Last,
$A^{(0,0)}_1 
= 
A^{(0)}_1 
=
\left(
\begin{array}{c}
0\\
A^{(0)}_{21}
\end{array}
\right)$,
by
Eq.~(\ref{A(0)}).
Substituting
these expressions for
$B^{f\perp}_0$
and
$A^{(0,0)}_1$
into
Eq.~(\ref{Q1s(0,0)}),
we obtain that
$Q_{1s}^{(0,0)} = 0$.

The induction hypothesis on $R_{(0)}$
can be verified
either by a similar argument,
or by recalling that
$R_{(0)} = (Q^{(0)})^{-1}$,
where $Q^{(0)}$
was shown above to be
block-triangular
to leading order.

\subsection{Proof of Theorem~\ref{t-CSPFq} for $q=1,2,\ldots$ \label{ss-q=1,2,...}}
We assume that
the induction hypothesis 
(\ref{Q(q)})--(\ref{R(q)}) 
holds for
$0,1,\ldots,q$
and
show that
it holds 
for $q+1$.
The proof proceeds
in four steps.
In step~1,
we derive explicit expressions for
$R_{(q+1)}$ and $Q^{(q+1)}$
in terms of
$R_{(q)}$ and $Q^{(q)}$;
these expressions also involve 
$U_{(q)}$
and $L_{(q)}$.
In step~2,
we derive 
the leading-order asymptotics
of $U_{(q)}$, and
in step~3
the leading-order asymptotics
of $L_{(q)}$.
Then, in step~4,
we substitute 
these results
into the expressions
derived in step~1 
to complete 
the induction.

\noindent \textbf{Step~1.}
We derive the expressions
for $Q^{(q+1)}$ and $R_{(q+1)}$.
Equations~(\ref{A-split})
and~(\ref{B-split}),
together with
the update formulas 
(\ref{A(q+1)})
for $A^{(q)}$
and~(\ref{B(q+1)})
for $B_{(q)}$,
yield 
\begin{eqnarray}
    Q^{(q+1)} 
&=&
    Q^{(q)} (I-\Uf_{(q)}) (I+\Lf_{(q)}) ,
\Label{Q-update}\\
    R_{(q+1)} 
&=&
    (I-\Lf_{(q)}) (I+\Uf_{(q)}) R_{(q)} .
\Label{R-update}
\end{eqnarray}
In terms of
the constituent blocks,
we have
\begin{eqnarray}
    Q_{1f}^{(q+1)}
&=&
      Q_{1f}^{(q)}
    + Q_{2f}^{(q)} L_{(q)} 
    - Q_{1f}^{(q)} U_{(q)} L_{(q)},
\Label{Q1f(q+1)}\\
    Q_{2f}^{(q+1)}
&=&
      Q_{2f}^{(q)} 
    - Q_{1f}^{(q)} U_{(q)} ,
\Label{Q2f(q+1)}\\
    Q_{1s}^{(q+1)}
&=&
      Q_{1s}^{(q)}
    + Q_{2s}^{(q)} L_{(q)}
    - Q_{1s}^{(q)} U_{(q)} L_{(q)} ,
\Label{Q1s(q+1)}\\
    Q_{2s}^{(q+1)}
&=&
      Q_{2s}^{(q)}
    - Q_{1s}^{(q)} U_{(q)} ,
\Label{Q2s(q+1)}
\end{eqnarray}
and
\begin{eqnarray}
    R^{1s\perp}_{(q+1)}
&=&
      R^{1s\perp}_{(q)}
    + U_{(q)} R^{2s\perp}_{(q)} ,
\Label{R1sperp(q+1)}\\
    R^{1f\perp}_{(q+1)}
&=&
      R^{1f\perp}_{(q)} 
    + U_{(q)} R^{2f\perp}_{(q)} ,
\Label{R1fperp(q+1)}\\
    R^{2s\perp}_{(q+1)}
&=&
      R^{2s\perp}_{(q)}
    - L_{(q)} R^{1s\perp}_{(q)}
    - L_{(q)} U_{(q)} R^{2s\perp}_{(q)} ,
\Label{R2sperp(q+1)}\\
    R^{2f\perp}_{(q+1)}
&=&
      R^{2f\perp}_{(q)}
    - L_{(q)} R^{1f\perp}_{(q)}
    - L_{(q)} U_{(q)} R^{2f\perp}_{(q)} .
\Label{R2fperp(q+1)}
\end{eqnarray}

\noindent \textbf{Step~2.}
We derive 
the leading-order asymptotics
of the matrix
$\Uf_{(q)}$.

Recall that 
$\Uf_{(q)}=
(\bfL^{11}_{(q)})^{-1}
\bfL^{12}_{(q)}$.
Moreover,
$\bfL^{11}_{(q)}$ 
is strictly
${\cal O}(1)$
and
$\bfL^{12}_{(q)}$
is strictly
${\cal O}(\eps^q)$
by Lemma~\ref{l-L-estim}.
Hence,
$\Uf_{(q)}
=
\Uf_{(q,q)} \eps^q
+ {\cal O}(\eps^{q+1})$,
with
$\Uf_{(q,q)} 
= 
(\bfL^{11}_{(q,0)})^{-1} \bfL^{12}_{(q,q)}$.
Therefore, 
it suffices
to derive 
the leading order asymptotics
of these blocks of $\bfL$. 

By definition,
$\bfL_{(q)} = B_{(q)} [A^{(q)},g]$.
Therefore,
\begin{eqnarray}
    \bfL_{(q)}
=
    \left(
          \begin{array}{cc}
            B^1_{(q)} [A_1^{(q)},g ]
          & B^1_{(q)} [A_2^{(q)},g ]\\
            B^2_{(q)} [A_1^{(q)},g ]
          & B^2_{(q)} [A_2^{(q)},g ]\\
          \end{array}
    \right) .
\Label{blox-L(q)}
\end{eqnarray}
The individual blocks of
$\bfL_{(q)}$
are obtained by
substituting Eqs.~(\ref{Ai-split}) 
and
(\ref{Bi-split}) into 
Eq.~(\ref{blox-L(q)}).
We observe that
one-half of all
the terms would vanish,
were they to
be evaluated on $\Me$,
by virtue of
Lemma~\ref{l-Ldiag}.
Since they are
evaluated on $\Kf^{(q+1)}$,
instead,
which is
${\mathcal O}(\eps^{q+1})$-close to $\Me$,
these terms are ${\mathcal O}(\eps^{q+2})$
and therefore of
higher order for
each of the blocks,
recall Lemma~\ref{l-L-estim}.
Thus,
\begin{eqnarray} \bfL^{11}_{(q)} 
&=& 
      R^{1s\perp}_{(q)} B^{s\perp} [A_f Q_{1f}^{(q)},g ] 
    + R^{1f\perp}_{(q)} B^{f\perp} [A_s Q_{1s}^{(q)},g ] , 
\Label{L11(q)-split}\\ 
\bfL^{12}_{(q)} 
&=& 
      R^{1s\perp}_{(q)} B^{s\perp} [A_f Q_{2f}^{(q)},g ]
    + R^{1f\perp}_{(q)} B^{f\perp} [A_s Q_{2s}^{(q)},g ] ,
\Label{L12(q)-split}\\
    \bfL^{21}_{(q)}
&=&
      R^{2s\perp}_{(q)} B^{s\perp} [A_f Q_{1f}^{(q)},g ]
    + R^{2f\perp}_{(q)} B^{f\perp} [A_s Q_{1s}^{(q)},g ] ,
\Label{L21(q)-split}
\end{eqnarray}
where the remainders of
${\cal O}(\eps^{q+2})$
have been ommited for brevity.
Recalling
the definition of
the Lie bracket,
we rewrite
Eq.~(\ref{L11(q)-split})
as
\begin{eqnarray}
    \bfL^{11}_{(q)} 
    &=&
    R^{1s\perp}_{(q)} B^{s\perp}
    \left((Dg) A_f Q_{1f}^{(q)}
          - \frac{d}{dt} \left(A_f Q_{1f}^{(q)}\right)
    \right) \nonumber\\
    &&\mbox{}
    + R^{1f\perp}_{(q)} B^{f\perp}
    \left((Dg) A_s Q_{1s}^{(q)}
          - \frac{d}{dt} \left(A_s Q_{1s}^{(q)}\right)
    \right),
\Label{L11(q)-split-exp}
\end{eqnarray}
where
we recall that
all of the quantities
are evaluated at
$(y,\psi_{(q+1)},\eps)$.
Next,
$(Dg) A_s$
and the two
time derivatives in
Eq.~(\ref{L11(q)-split-exp})
are zero to leading order
by Lemma~\ref{l-DgoAo}
and
\cite[Lemma A.2]{ZKK-2003},
respectively.
Therefore,
to leading order
Eq.~(\ref{L11(q)-split-exp})
becomes
\begin{eqnarray}
    \bfL^{11}_{(q,0)} 
&=& 
    R^{1s\perp}_{(q,0)} 
    B^{s\perp}_0
    (Dg)_0 
    A_f^0 
    Q_{1f}^{(q,0)}.
\Label{L11(q)-split-lead}
\end{eqnarray}
Here,
$\bfL^{11}_{(q,0)}$
stands for 
the leading-order
term
in the 
asymptotic expansion of
$\bfL^{11}_{(q)}(y,\psi_{(q+1)}(y),\eps)$,
and
the right member is
the leading order term in
the asymptotic expansion of
$(R^{1s\perp}_{(q)} B^{s\perp} (Dg) A_f Q_{1f}^{(q)}) (y,h_\eps(y),\eps)$.

We derive
a similar formula for
$\bfL^{12}_{(q,q)}$.
First,
we rewrite
Eq.~(\ref{L12(q)-split})
as
\begin{eqnarray}
    \bfL^{12}_{(q)} 
&=&
    R^{1s\perp}_{(q)} B^{s\perp}
    \left((Dg) A_f Q_{2f}^{(q)}
           - \frac{d}{dt} \left(A_f Q_{2f}^{(q)}\right)
    \right) 
\nonumber\\
&&\mbox{} 
    + R^{1f\perp}_{(q)} B^{f\perp}
    \left((Dg) A_s Q_{2s}^{(q)}
           - \frac{d}{dt} \left(A_s Q_{2s}^{(q)}\right)
    \right) .
\Label{L12(q)-split-exp}
\end{eqnarray}
Next,
$Q_{2f}^{(q)} = {\cal O}(\eps^q)$,
$Q_{2s}^{(q)} = {\cal O}(1)$,
$R_{(q)}^{1s\perp} = {\cal O}(1)$,
and
$R_{(q)}^{1f\perp} = {\cal O}(\eps^q)$,
by the induction hypothesis
(\ref{Q(q)})--(\ref{R(q)}).
Thus,
\cite[Lemma A.2]{ZKK-2003}
implies that
the two terms
in Eq.~(\ref{L12(q)-split-exp})
involving 
time derivatives
are ${\cal O}(\eps^{q+1})$
and therefore of higher order.
Also,
$(Dg) A_s$ is
zero to leading order
by Lemma~\ref{l-DgoAo},
and thus
\begin{eqnarray}
    \bfL^{12}_{(q,q)} 
    &=&
    R_{(q,0)}^{1s\perp}
    B^{s\perp}_{0}
    (Dg)_0
    A_{f}^0
    Q^{(q,q)}_{2f} .
\Label{L12(q)-split-lead}
\end{eqnarray}

We now
substitute
$\bfL^{11}_{(q,0)}$ 
and $\bfL^{12}_{(q,q)}$
from
Eqs.~(\ref{L11(q)-split-lead})
and
(\ref{L12(q)-split-lead})
in the expression
$U_{(q,q)}
=
(\bfL^{11}_{(q,0)})^{-1} 
\bfL^{12}_{(q,q)} $
to find
the desired expression 
for
$U_{(q,q)}$
in terms of
$Q^{(q)}$,
\be
    U_{(q,q)}
=  
    \left(Q^{(q,0)}_{1f}\right)^{-1} Q^{(q,q)}_{2f} .
\Label{U(q,q)-Q}
\ee

We also need 
an expression for
$U_{(q,q)}$
in terms of
blocks of $R_{(q)}$,
which we will use in
Eqs.~(\ref{R1sperp(q+1)})--(\ref{R2fperp(q+1)}).
Since
$R_{(q)}$ has
the near
block-diagonal structure
given by 
the induction hypothesis
(\ref{Q(q)})--(\ref{R(q)})
and
$Q^{(q)}$ is
its inverse,
we find
\begin{eqnarray}
    \hspace{-2.0em}Q^{(q)}
&\hspace{-0.75em}=\hspace{-0.75em}& 
    \left(
           \begin{array}{cc}
             (R_{(q,0)}^{1s\perp})^{-1} 
           & - \eps^q 
               (R_{(q,0)}^{1s\perp})^{-1}
               R_{(q,q)}^{1f\perp}
               (R_{(q,0)}^{2f\perp})^{-1} \\
             - \eps^{q+1} 
               (R_{(q,0)}^{2f\perp})^{-1}
               R_{(q,q+1)}^{2s\perp}
               (R_{(q,0)}^{1s\perp})^{-1}
           & (R_{(q,0)}^{2f\perp})^{-1}
           \end{array} 
    \right) \hspace{-0.25em},~~~
\Label{Q(q)-vs-R(q)}
\end{eqnarray}
to leading order
for each of the blocks
and
for $q = 1,2,\ldots\,$.
Equations~(\ref{U(q,q)-Q})
and~(\ref{Q(q)-vs-R(q)})
lead to
the desired expression 
for
$U_{(q,q)}$
in terms of
$R_{(q)}$,
\be
    U_{(q,q)}
=  
    - R_{(q,q)}^{1f\perp}
      \left(R_{(q,0)}^{2f\perp}\right)^{-1} .
\Label{U(q,q)-R}
\ee

\noindent
\textbf{Step~3.}
We derive the leading-order asymptotics
of the matrix
$L_{(q)}$.

Recall that 
$\Lf_{(q)}=
\bfL^{21}_{(q)}
(\bfL^{11}_{(q)})^{-1}$.
Moreover,
by Lemma~\ref{l-L-estim},
$\bfL^{11}_{(q)}$ 
is strictly
${\cal O}(1)$
and
$\bfL^{21}_{(q)}$
is strictly
${\cal O}(\eps^{q+1})$.
Hence,
$\Lf_{(q)}
=
\Lf_{(q,q+1)} \eps^{q+1}
+ {\cal O}(\eps^{q+2})$,
with
$\Lf_{(q,q+1)} 
= 
\bfL^{21}_{(q,q+1)}
\bfL^{11}_{(q,0)})^{-1}. $
An expression for
$\bfL^{11}_{(q,0)}$
was derived in 
Eq.~(\ref{L11(q)-split-lead}),
so here
we focus on
$\bfL^{21}_{(q,q+1)}$.

Equation~(\ref{L21(q)-split})
and
the definition of
the Lie bracket
imply that 
\begin{eqnarray}
    \bfL^{21}_{(q)}
&=& 
    R_{(q)}^{2s\perp} B^{s\perp} 
    \left(                        
          (Dg) A_f Q^{(q)}_{1f}
          - \frac{d}{dt} \left(A_f Q^{(q)}_{1f}\right)
    \right) \nonumber\\
    &&\mbox{}
    + R_{(q)}^{2f\perp} B^{f\perp} 
    \left(
          (Dg) A_s Q^{(q)}_{1s}
          - \frac{d}{dt} \left(A_s Q^{(q)}_{1s}\right)
    \right) .
\Label{L21(q)-split-exp}
\end{eqnarray}
Next,
$Q^{(q)}_{1f} = {\cal O}(1)$,
$Q^{(q)}_{1s} = {\cal O}(\eps^{q+1})$,
$R_{(q)}^{2s\perp} = {\cal O}(\eps^{q+1})$,
and
$R_{(q)}^{2f\perp} = {\cal O}(1)$,
by the induction hypothesis.
Also,
the time derivatives are
${\cal O}(\eps)$ by
\cite[Lemma A.2]{ZKK-2003},
and thus
the two terms
in Eq.~(\ref{L21(q)-split-exp})
that involve
time derivatives
are ${\cal O}(\eps^{q+2})$.
Last,
$(Dg) A_s = {\cal O}(\eps)$
by Lemma~\ref{l-DgoAo}.
Thus,
we find
\begin{eqnarray}
    \bfL^{21}_{(q,q+1)} 
    &=&
    R^{2s\perp}_{(q,q+1)}
    B^{s\perp}_{0}
    (Dg)_0
    A_{f}^0
    Q_{1f}^{(q,0)}.
\Label{L21(q)-split-lead}
\end{eqnarray}
Equations~(\ref{L11(q)-split-lead})
and
(\ref{L21(q)-split-lead})
yield
the desired formula for
$L_{(q,q+1)}$
in terms of
the blocks of
$R_{(q)}$,
\be
    L_{(q,q+1)}
=
    \bfL^{21}_{(q,q+1)} \left(\bfL^{11}_{(q,0)}\right)^{-1} 
=
    R_{(q,q+1)}^{2s\perp}
    \left(R_{(q,0)}^{1s\perp}\right)^{-1} .
\Label{L(q,q+1)-R}
\ee

Next,
we recast
Eq.~(\ref{L(q,q+1)-R})
in terms of
blocks of $Q^{(q)}$,
in order to use it in
Eqs.~(\ref{Q1f(q+1)})--(\ref{Q2s(q+1)}).
The matrix
$R_{(q)}$ is
the inverse of
$Q^{(q)}$
and
has
the near 
block-diagonal form
given
in~(\ref{R(q)}).
Thus,
\begin{eqnarray}
    \hspace{-2.0em}R_{(q)}
&\hspace{-0.75em}=\hspace{-0.75em}& 
    \left(
           \begin{array}{cc}
             (Q^{(q,0)}_{1f})^{-1} 
           & - \eps^q
               (Q^{(q,0)}_{1f})^{-1}
               Q^{(q,q)}_{2f}
               (Q^{(q,0)}_{2s})^{-1} \\
             - \eps^{q+1}
               (Q^{(q,0)}_{2s})^{-1}
               Q^{(q,q+1)}_{1s}
               (Q^{(q,0)}_{1f})^{-1}
           & (Q^{(q,0)}_{2s})^{-1}
           \end{array} 
    \right) \hspace{-0.3em},~~~~
\Label{R(q)-vs-Q(q)}
\end{eqnarray}
to leading order
for each block
and
for $q = 1,2,\ldots\,$.
Equations~(\ref{L(q,q+1)-R})
and~(\ref{R(q)-vs-Q(q)})
lead to
the desired expression
for
$L_{(q,q+1)}$
in terms of
the blocks of
$Q^{(q)}$,
\be
    L_{(q,q+1)}
=
    - \left(Q^{(q,0)}_{2s}\right)^{-1}
      Q^{(q,q+1)}_{1s} .
\Label{L(q,q+1)-Q}
\ee

\noindent
\textbf{Step 4.}
We substitute 
the results 
obtained in Step~2
and Step~3 
into the formulas
(\ref{Q1f(q+1)})--(\ref{R2fperp(q+1)})
derived in Step~1.

Equations~(\ref{Q1f(q+1)})
and
(\ref{Q2s(q+1)}),
together with
the induction hypothesis
and
the estimates
$U_{(q)} = {\cal O}(\eps^q)$ and
$L_{(q)} = {\cal O}(\eps^{q+1})$,
imply that
$Q_{1f}^{(q+1)}$
and
$Q_{2s}^{(q+1)}$
remain ${\cal O}(1)$.
This concludes
the estimation of
these blocks.

Next,
we show that
$Q_{2f}^{(q+1)} =
{\cal O}(\eps^{q+1})$.
First,
$Q_{2f}^{(q+1)}$
and
$Q_{2f}^{(q)}$
are equal
up to and including
terms of ${\cal O}(\eps^{q-1})$,
by 
Eq.~(\ref{Q2f(q+1)})
and
the estimate on
$U_{(q)}$.
Thus,
$Q_{2f}^{(q+1,i)} = 0$
for $i = 0,1,\ldots\,,q-1$,
by the induction hypothesis
on $Q_{2f}^{(q)}$.
It suffices to show that
$Q_{2f}^{(q+1,q)} = 0$.
Equation~(\ref{Q2f(q+1)})
implies that
\be
    Q_{2f}^{(q+1,q)}
=
      Q_{2f}^{(q,q)}
    - Q_{1f}^{(q,0)} U_{(q,q)} .
\Label{Q2f(q+1,q)}
\ee
The right member of
this equation
is zero,
by Eq.~(\ref{U(q,q)-Q}),
and 
the estimation of
$Q_{2f}^{(q+1)}$
is complete.

Finally,
we show that 
$Q_{1s}^{(q+1)} =
{\cal O}(\eps^{q+2})$
to complete the estimates on
the blocks of $Q^{(q+1)}$.
First,
$Q_{1s}^{(q+1)}$
and
$Q_{1s}^{(q)}$
are equal
up to and including
terms of ${\cal O}(\eps^q)$,
by 
Eq.~(\ref{Q1s(q+1)})
and
the order estimates on
$U_{(q)}$ and $L_{(q)}$.
Thus,
$Q_{1s}^{(q+1,i)} = 0$
for $i = 0,1,\ldots,q$,
by the induction hypothesis
on $Q_{1s}^{(q)}$.
It suffices to
show that
$Q_{1s}^{(q+1,q+1)} = 0$.
Equation~(\ref{Q1s(q+1)})
implies that
\be
    Q_{1s}^{(q+1,q+1)}
=
      Q_{1s}^{(q,q+1)}
    + Q_{2s}^{(q,0)} L_{(q,q+1)} ,
\Label{Q1s(q+1,q+1)}
\ee
where
the right member of
this equation
is zero
by Eq.~(\ref{L(q,q+1)-Q}).
The estimation of
$Q_{1s}^{(q+1)}$
is complete.

The blocks of
$R_{(q)}$
may be estimated
in an entirely similar manner,
using
Eqs.~(\ref{R1sperp(q+1)})--(\ref{R2fperp(q+1)}),
instead of
Eqs.~(\ref{Q1f(q+1)})--(\ref{Q2s(q+1)}),
and
Eqs.~(\ref{U(q,q)-R}) and (\ref{L(q,q+1)-R}),
instead of
Eqs.~(\ref{U(q,q)-Q}) and (\ref{L(q,q+1)-Q}).
The proof of Theorem~\ref{t-CSPFq} is complete.

\section{The Michaelis--Menten--Henri Model \label{s-MMH model}}
\setcounter{equation}{0}
In this section, 
we
illustrate
Theorem~\ref{t-CSPFq}
by applying
the CSP method to 
the Michaelis--Menten--Henri 
(MMH) mechanism of 
enzyme kinetics 
\cite{P-1987, PL-1984}.
We consider
the planar system of 
ODEs for 
a slow variable $s$ 
and
a fast variable $c$,
\begin{eqnarray}
    s' &=& \eps (- s + (s + \kappa - \lambda)c) , \Label{MMH-s} \\
    c' &=& s - (s + \kappa)c . \Label{MMH-c}
\end{eqnarray}
The parameters
satisfy the inequalities
$0 < \eps \ll 1$ and $\kappa > \lambda > 0$.
Only nonnegative values of 
$s$ and $c$ 
are relevant.
The system of 
Eqs.~(\ref{MMH-s})--(\ref{MMH-c})
is of 
the form~(\ref{eq-y})--(\ref{eq-z})
with 
$m=1$, $n=1$,
$y=s$, $z=c$,
$g_1 = -s+(s+\kappa-\lambda)c$,
and
$g_2 = s-(s+\kappa)c$.

\subsection{Slow Manifolds and Fast Fibers \label{ss-MMH-Me&Fe}}
In the limit as $\eps \downarrow 0$,
the dynamics of the MMH equations
are confined to the reduced slow manifold
\be
    \Mo = \{ (c, s) : c = \frac{s}{s+\kappa} , s \ge 0 \} . \quad
\Label{MMH-Mo}
\ee
The manifold $\Mo$ is asymptotically stable,
so
there exists
a locally invariant
slow manifold $\Me$
for all sufficiently small $\eps$
that is $\mathcal{O}(\eps)$ close
to $\Mo$ on any compact set.
Moreover, $\Me$ is the graph
of a function $h_\eps$,
\be
    \Me = \{ (c, s) : c = h_\eps (s) , s \ge 0 \} ,
    \Label{MMH-Me}
\ee
and $h_\eps$ admits an asymptotic expansion,
$h_\eps = h_0 + \eps h_1 + \eps^2 h_2 + \cdots \,$.
The coefficients are found from
the invariance equation,
\be
  s - (s + \kappa) h_\eps (s)
  =
  \eps h_\eps' (s) (-s + (s + \kappa - \lambda) h_\eps (s)) .
\ee
The first few coefficients are
\be
    h_0 (s) = \frac{s}{s+\kappa} , \quad
    h_1(s) = \frac{\kappa\lambda s}{(s+\kappa)^4} , \quad
    h_2(s) = \frac{\kappa\lambda s(2\kappa\lambda-3\lambda s-\kappa s-\kappa^2)}{(s+\kappa)^7} .
\Label{Me-MMH}
\ee

In the limit as $\eps \downarrow 0$,
each 
line of constant $s$
is trivially invariant 
under 
Eqs.~(\ref{MMH-s})--(\ref{MMH-c}).
These are the 
(one-dimensional) fast fibers
$\Fo^p$
with base point 
$p = (s,h_0(s))\in\Mo$.
All points on $\Fo^p$ 
contract exponentially fast to $p$
with rate constant
$-(s+\kappa)$.
The fast fiber $\Fo^p$
perturbs to a curve $\Fe^p$
that is
${\cal O}(\eps)$ close to 
$\Fo^p$
in any compact neighborhood
of $\Me$.
The fast fibers $\Fe^p$,
$p \in \Me$,
form an invariant family.

\subsection{Asymptotic Expansions of the Fast Fibers}
To derive asymptotic information 
about the fast fibers,
we look for
general solutions of 
Eqs.~(\ref{MMH-s})--(\ref{MMH-c})
that are given by
asymptotic expansions,
\be
    s(t;\eps) = \sum_{i=0} \eps^i s_i(t) , \quad
    c(t;\eps) = \sum_{i=0} \eps^i c_i(t),
\Label{sc-exp} 
\ee
where the coefficients 
$s_i$ and $c_i$ 
are determined 
order by order. 

Consider the fast fiber $\Fe^p$
with base point $p=(s,h_\eps(s))$,
and let 
$(s^A,c^A)$
and
$(s^B,c^B)$ 
be two points on it; 
let $\Delta s(t) = s^B(t)-s^A(t)$ 
and $\Delta c(t) = c^B(t)-c^A(t)$. 
The distance 
between any two points 
on the same fast fiber 
will contract exponentially fast 
towards zero
at the
${\cal O}(1)$ rate,
as long as 
both points are chosen
in a neighborhood of
$\Me$. 
We may write
\be
    \Delta s(t;\eps) = \sum_{i=0} \eps^i \Delta s_i(t) , \quad
    \Delta c(t;\eps) = \sum_{i=0} \eps^i \Delta c_i(t),
\Label{DsDc-exp} 
\ee
where
$\Delta s_i(t) = s_i^B(t)-s_i^A(t)$
and
$\Delta c_i(t) = c_i^B(t)-c_i^A(t)$. The condition on
fast exponential decay of $\Delta s(t)$ and $\Delta c(t)$ translates into
\be
  \Delta s_i(t) = {\cal O}(e^{-C_s t}) , \;
  \Delta c_i(t) = {\cal O}(e^{-C_c t}) , \quad
 t\to\infty ,
\Label{genl sici-cond}
\ee
for some positive constants 
$C_s$ and $C_c$.
We let
$(s^A,c^A)$ and $(s^B,c^B)$ be 
infinitesimally close, 
since we are interested in vectors 
tangent to the fast fiber.

\subsubsection{${\cal O}(1)$ Fast Fibers}
Substituting the expansions
(\ref{sc-exp})
into Eqs.~(\ref{MMH-s})--(\ref{MMH-c})
and equating $\mathcal{O} (1)$ terms,
we find 
\begin{eqnarray}
    s_0' &=& 0 ,
\Label{s0-ODE} \\
    c_0' &=& s_0-(s_0+\kappa) c_0. 
\Label{c0-ODE} 
\end{eqnarray}
The equations can be integrated,
\begin{eqnarray}
    s_0(t) &=& s_0(0)=s_0 , 
\Label{s0-soln}\\
    c_0(t) &=& \frac{s_0}{s_0+\kappa}+
  \left(
   c_0(0)-\frac{s_0}{s_0+\kappa} 
  \right) 
  e^{-(s_0+\kappa)t} . 
\Label{c0-soln} 
\end{eqnarray}
Hence,
\begin{eqnarray}
      \Delta s_0(t) &=& \Delta s_0(0) , 
\Label{Ds0-soln} \\
      \Delta c_0(t) &=& 
      \Delta c_0(0) e^{-(s_0+\kappa)t} 
    + \left(\partial_{s_0} c_0(t)\right) 
      \Delta s_0(0)
    + {\cal O}((\Delta s_0(0))^2) . 
\Label{Dc0-soln} 
\end{eqnarray}
The points 
$A$ and $B$ 
lie on 
the same fiber
if and only if
\be
    \Delta s_0(0) = 0. \Label{cond0}
\ee
Thus,
Eq.~(\ref{Dc0-soln}) simplifies to
\be
    \Delta c_0(t) = \Delta c_0(0) e^{-(s_0+\kappa)t} ,
\ee
and
$\Delta c_0(t)$ decays 
exponentially towards zero, 
irrespective of 
the choice of $\Delta c_0(0)$.
Hence, 
$\Delta c_0(0)$ is 
a free parameter.

We conclude that,
to $\mathcal{O} (1)$,
any vector
$\left( \begin{array}{c} 0 \\ \alpha \end{array} \right)$
with $\alpha$ constant ($\alpha\neq 0$)
is tangent to every fast fiber
at the base point.

\subsubsection{${\cal O}(\eps)$ Fast Fibers}
At ${\cal O}(\eps)$, 
we obtain the equations
\begin{eqnarray}
    s_1' &=& -s_0+(s_0+\kappa-\lambda)c_0 , 
\Label{s1-ODE} \\ 
    c_1' &=& s_1-(s_0+\kappa)c_1-s_1c_0. 
\Label{c1-ODE}
\end{eqnarray}
Using Eqs.~(\ref{s0-soln}) and~(\ref{c0-soln}),
we integrate Eq.~(\ref{s1-ODE})
to obtain
\be
s_1(t) = s_1(0)-\frac{\lambda s_0}{s_0+\kappa}t +
       \frac{s_0+\kappa-\lambda}{s_0+\kappa}
       \left(
           c_0(0)-\frac{s_0}{s_0+\kappa} 
       \right)
       (1-e^{-(s_0+\kappa)t}). \Label{s1-soln}
\ee
Therefore, 
at ${\cal O}(\eps)$, 
\be
    \Delta s_1(t) = \Delta s_1(0)+
                \frac{s_0+\kappa-\lambda}{s_0+\kappa}
                \Delta c_0(0) (1-e^{-(s_0+\kappa)t}) .
\ee
For 
the two points to
have the same phase 
asymptotically,
it is 
necessary that
$\lim_{t\to\infty} \Delta s_1(t) = 0$.
This condition
is satisfied if and only if
\be
    \Delta s_1(0) = -\frac{s_0+\kappa-\lambda}{s_0+\kappa} \Delta c_0(0). \Label{cond1}
\ee
Next,
$c_1(t)$ follows
upon integration of Eq.~(\ref{c1-ODE}),
\begin{eqnarray}
c_1(t) &=& c_1(0)e^{-(s_0+\kappa)t}\nonumber\\
       &&\mbox{}
       +\frac{\kappa}{(s_0+\kappa)^2}
       \left(
           s_1(0) +
           \frac{s_0+\kappa-\lambda}{s_0+\kappa}
           \left(
               c_0(0)-\frac{s_0}{s_0+\kappa} 
           \right)
       \right)
       (1-e^{-(s_0+\kappa)t})\nonumber\\
       &&\mbox{}
        -\left(
              c_0(0)-\frac{s_0}{s_0+\kappa} 
        \right)
        \left(
          s_1(0) +
          \frac{s_0+\kappa-\lambda}{s_0+\kappa}
          \left(
            c_0(0)+\frac{\kappa-s_0}{s_0+\kappa} 
          \right)
        \right)
       te^{-(s_0+\kappa)t}\nonumber\\
       &&\mbox{}
        -\frac{s_0+\kappa-\lambda}{(s_0+\kappa)^2}
        \left(
          c_0(0)-\frac{s_0}{s_0+\kappa} 
        \right)^2
       (e^{-2(s_0+\kappa)t}- e^{-(s_0+\kappa)t})\nonumber\\
       &&\mbox{}
        +\frac{\lambda s_0}{2(s_0+\kappa)}
        \left(
          c_0(0)-\frac{s_0}{s_0+\kappa} 
        \right)
       t^2 e^{-(s_0+\kappa)t}\nonumber\\
       &&\mbox{}
       -\frac{\kappa \lambda s_0}{(s_0+\kappa)^4}
       (e^{-(s_0+\kappa)t}+(s_0+\kappa) t-1).~
\Label{c1-soln}
\end{eqnarray}
We infer from 
this expression that
$\lim_{t\to\infty}\Delta c_1(t)=0$,
as long as 
Eqs.~(\ref{cond1}) 
and (\ref{cond0})
hold. 
Hence, 
$\Delta c_1(0)$ is a free parameter, 
just like $\Delta c_0(0)$, 
and the only condition 
that arises at ${\cal O}(\eps)$
is (\ref{cond1}) 
on $\Delta s_1(0)$.

We conclude that
any vector
\begin{eqnarray}
\left(
    \begin{array}{c}
    0\\
    \alpha
    \end{array}
\right)
&+&
\eps\left(
       \begin{array}{c}
       - \left(1-\frac{\lambda}{s_0+\kappa}\right) \alpha \\
       \beta
       \end{array}
   \right) ,
\Label{MMH-TF1}
\end{eqnarray}
with $\alpha$ and $\beta$ constant
($\alpha \neq 0$),
is tangent to every fast fiber
at the base point
up to and including terms
of $\mathcal{O} (\eps)$.
Any such vector
may be written as
the product of 
a free parameter
and
a constant vector
(fixed by $s_0$),
\begin{eqnarray}
(\alpha + \eps \beta)
\left(
    \begin{array}{c}
       - \eps \left(1-\frac{\lambda}{s_0+\kappa}\right) \\
    1
    \end{array}
\right)
+ {\mathcal O}(\eps^2).
\end{eqnarray}

\subsubsection{${\cal O}(\eps^2)$ Fast Fibers}
At {${\cal O}(\eps^2)$, we obtain the equation
\be 
    s_2' = s_1 (c_0-1)+(s_0+\kappa-\lambda)c_1. \Label{s2-ODE} 
\ee
Direct integration yields
\begin{eqnarray}
    s_2(t) 
    &=& 
%
%
      s_2(0)
+ \left[
        \frac{\lambda}{(s_0+\kappa)^2}
        \left(
              c_0(0)-\frac{s_0}{s_0+\kappa} 
        \right)
        - \frac{\kappa(s_0+\kappa-\lambda)}{(s_0+\kappa)^3}
\right] s_1(0)
\nonumber\\
   &&\mbox{}
- \left[
      \frac{\kappa(s_0+\kappa-\lambda)(s_0+\kappa-2\lambda) + \lambda^2 s_0}
           {(s_0+\kappa)^4}
\right]
      \left(
            c_0(0)-\frac{s_0}{s_0+\kappa} 
      \right)
\nonumber\\
    &&\mbox{}
    + \frac{\lambda(s_0+\kappa-\lambda)}{2(s_0+\kappa)^3}  
      \left(
            c_0(0)-\frac{s_0}{s_0+\kappa} 
      \right)^2
\nonumber\\
       &&\mbox{}
    + \left(1-\frac{\lambda}{s_0+\kappa}\right)
      \left(c_1(0) - \frac{\kappa\lambda s_0}{(s_0+\kappa)^4}\right)
\nonumber\\
       &&\mbox{}
%
%
- \frac{\kappa \lambda}{(s_0+\kappa)^2}
\Bigg[
      s_1(0)
    + \frac{s_0+\kappa-\lambda}{s_0+\kappa} 
      \left(
            c_0(0)-\frac{2 s_0}{s_0+\kappa} 
      \right)
\Bigg] t
\nonumber\\
       &&\mbox{}
%
%
    + \frac{\kappa \lambda^2 s_0}{2(s_0+\kappa)^3} t^2
    + {\cal R}(t) ,
\Label{s2-soln}
\end{eqnarray}
where the remainder ${\cal R}(t)$
involves the functions
$e^{-(s_0+\kappa) t}$,
$t e^{-(s_0+\kappa) t}$,
$t^2 e^{-(s_0+\kappa) t}$,
and
$e^{-2 (s_0+\kappa) t}$. 
From this expression we find
\begin{eqnarray}
    \Delta s_2(t)
&=&  
      \Delta s_2(0)
    + \left(\partial_{s_0} s_2(t)\right)
      \Delta s_0(0)
    + \left(\partial_{c_0} s_2(t)\right)
      \Delta c_0(0)
\nonumber\\
    &&\mbox{}
    + \left(\partial_{s_1} s_2(t)\right)
      \Delta s_1(0)
    + \left(\partial_{c_1} s_2(t)\right)
      \Delta c_1(0)
    + {\cal O}(2) 
    + {\cal O}(e^{-Ct}) ,~~~~~
\Label{Ds2-Taylor}
\end{eqnarray}
for some
$C>0$.
Here,
$\partial_{c_0}$
is an abbreviation for
the partial derivative
$\partial_{c_0(0)}$,
and so on,
and
${\cal O}(2)$ denotes 
quadratic terms in
the multivariable
Taylor expansion.
First,
we recall that
$\Delta s_0(0) = 0$
by Eq.~(\ref{cond0}).
Next,
we calculate
the partial derivatives in
each of the
three remaining terms,
\begin{eqnarray}
    \partial_{c_0} s_2(t)
&=&
    \frac{\lambda s_1(0)}{(s_0+\kappa)^2}
    - \frac{\kappa(s_0+\kappa-\lambda)(s_0+\kappa-2\lambda)
            + \lambda^2 s_0}
           {(s_0+\kappa)^4}
\nonumber\\
    &&\mbox{}
    + \frac{\lambda(s_0+\kappa-\lambda)}{(s_0+\kappa)^3}
    \left(
          c_0(0)-\frac{s_0}{s_0+\kappa} 
    \right)
    - \frac{\kappa\lambda(s_0+\kappa-\lambda)}{(s_0+\kappa)^3} t ,~~~
\Label{Dc0}\\
    \partial_{s_1} s_2(t)
&=&
    \frac{\lambda}{(s_0+\kappa)^2}
    \left(
          c_0(0)-\frac{s_0}{s_0+\kappa} 
    \right)
    - \frac{\kappa(s_0+\kappa-\lambda)}{(s_0+\kappa)^3}
\nonumber\\
    &&\mbox{}
    - \frac{\kappa\lambda}{(s_0+\kappa)^2} t ,
\Label{Ds1}\\
    \partial_{c_1} s_2(t)
&=&
    1 - \frac{\lambda}{s_0+\kappa} .
\Label{Dc1}
\end{eqnarray}
We substitute
these expressions
into
Eq.~(\ref{Ds2-Taylor}),
recall Eq.~(\ref{cond1}),
and carry out the algebra
to obtain
\begin{eqnarray}
    \Delta s_2(t)
&=&  
      \Delta s_2(0)
    + \left(1- \frac{\lambda}{s_0+\kappa}\right) 
      \Delta c_1(0)
\nonumber\\
&&\mbox{}
    + \frac{\lambda}{(s_0+\kappa)^2}
      \left(
              s_1(0)
            + \frac{\kappa(s_0+\kappa-\lambda)-\lambda s_0}
                   {(s_0+\kappa)^2} 
      \right)
      \Delta c_0(0)
\nonumber\\
&&\mbox{}
      + {\cal O}(2)
      + {\cal O}(e^{-Ct}), \quad C>0.
\Label{Ds2-soln}
\end{eqnarray}

In the limit $t\to\infty$,
Eq.~(\ref{Ds2-soln})
yields the condition 
\begin{eqnarray}
  \Delta s_2(0)
  &=& 
  - \left(1-\frac{\lambda}{s_0+\kappa}\right) \Delta c_1(0) \nonumber\\
  &&\mbox{}
  -\frac{\lambda}{(s_0+\kappa)^2}\left(
  s_1(0) + \frac{\kappa(s_0+\kappa-\lambda)-\lambda s_0}{(s_0+\kappa)^2} 
                                \right)\Delta c_0(0).
\Label{cond2}
\end{eqnarray}
Finally, 
$\Delta c_2(t)$ 
vanishes exponentially,
as follows directly
from the conditions 
(\ref{cond1}) and (\ref{cond2}).
Thus, 
no further conditions
besides (\ref{cond2})
arise at ${\cal O}(\eps^2)$.

We conclude that
any vector 
\begin{eqnarray}
\left(
    \begin{array}{c}
    0\\
    \alpha
    \end{array}
\right)
&+&
\eps\left(
       \begin{array}{c}
       - \left(1-\frac{\lambda}{s_0+\kappa}\right) \alpha \\
       \beta
       \end{array}
   \right)\nonumber\\
&+&
\eps^2\left(
         \begin{array}{c}
         - \left(1-\frac{\lambda}{s_0+\kappa}\right) \beta
         -\frac{\lambda}{(s_0+\kappa)^2}\left(
                                s_1(0) + \frac{\kappa(s_0+\kappa-\lambda)-\lambda s_0}{(s_0+\kappa)^2}
                            \right)\alpha\\
         \gamma
         \end{array}
    \right) ,~~~~~~~
\Label{MMH-TF2}
\end{eqnarray}
with 
$\alpha$, $\beta$, and $\gamma$ 
constant
($\alpha\neq0$), 
is tangent to 
every fiber at the base point,
up to and including terms
of ${\cal O}(\eps^2)$.

\subsection{CSP Approximations of the Fast Fibers \label{ss-MMH}}
We choose the stoichiometric vectors
as the basis vectors, so
\be
    \Af^{(0)} 
    = 
    (\Af_1^{(0)}, \, \Af_2^{(0)})
    = 
    \left(
          \begin{array}{cc}
          0 & 1 \\
          1 & 0
          \end{array}
    \right) , \quad
    \Bf_{(0)}
    =
    \left(
          \begin{array}{cc}
          \Bf^1_{(0)} \\
          \Bf^2_{(0)}
          \end{array}
     \right)
     =
     \left(
          \begin{array}{cc}
          0 & 1 \\
          1 & 0
          \end{array}
      \right) .
\Label{1/A,B}
\ee
The CSP condition
$B^1_{(0)} g = 0$
is satisfied if
$c = h_0 (s)$,
so the CSP manifold
$\Kf^{(0)}$
coincides with $\Mo$.
With this choice of
initial basis, we have
\be
    \bfL_{(0)}
    = 
    \Bf_{(0)} (Dg) \Af^{(0)} 
    = 
    \left(
          \begin{array}{cc}
            - (s+\kappa) 
          & - (c-1) \\
            \eps(s+\kappa-\lambda) 
          & \eps(c-1)
          \end{array}
    \right) .
\Label{MMH-Lambda(0)}
\ee

\subsubsection{First Iteration}
At any point $(s,c)$,
we have
\begin{eqnarray}
    \Af_1^{(1)}
&=&
    \left(
          \begin{array}{c}
          0\\
          1
          \end{array}
    \right)
    +
    \eps
    \frac{s+\kappa-\lambda}{s+\kappa}
    \left(
          \begin{array}{c}
           -1\\
          \mbox{} \frac{c-1}{s+\kappa}
          \end{array}
    \right) , \qquad
    \Af_2^{(1)} 
=
    \left(
          \begin{array}{c}
            1 \\
          \mbox{}- \frac{c-1}{s+\kappa}
          \end{array}
    \right) ,
\Label{A(1)-genl}\\
    \Bf^1_{(1)} 
&=&
    \left( -\Af_{22}^{(1)} , \, \Af_{12}^{(1)} \right) , 
\hskip1.45truein
    \Bf^2_{(1)} 
=
    \left( \Af_{21}^{(1)} , \, -\Af_{11}^{(1)} \right) .~~~
\Label{B(1)-genl}
\end{eqnarray}
In the first step,
we evaluate
$\Af_2^{(1)}$
and
$\Bf^1_{(1)}$
on $\Kf^{(0)}$
to obtain
\be
    \Af_2^{(1)} 
=
    \left(
          \begin{array}{c}
          1\\
          \frac{\kappa}{(s+\kappa)^2}
          \end{array}
    \right),
\quad
    \Bf^1_{(1)} 
=
    \left(
          \mbox{}-\frac{\kappa}{(s+\kappa)^2}, \, 1
    \right) .
\Label{A2(1)-B1(1)}
\ee
Hence,
the CSP condition,
\be
    \Bf^1_{(1)} g
=  
      s
    - (s+\kappa)c
    - \eps\frac{\kappa(-s+(s+\kappa-\lambda)c)}{(s+\kappa)^2}
=
    0 ,
\Label{f1(1)=0}
\ee
is satisfied if
\be
    c 
= 
      \frac{s}{s+\kappa}
    + \eps
      \frac{\kappa\lambda s}
           {(s+\kappa)^4}
    - \eps^2
      \frac{\kappa^2\lambda s(s+\kappa-\lambda)}
           {(s+\kappa)^7}
    + \mathcal{O}(\eps^3) .
\Label{proj2}
\ee
Equation~(\ref{proj2}) 
defines $\Kf^{(1)}$,
the CSPM of order one,
which agrees with $\Me$
up to and including terms 
of $\mathcal{O}(\eps)$;
recall Eq.~(\ref{Me-MMH}).

Then,
in the second step,
the new fast basis vector, 
$A_1^{(1)}$, 
and its complement,
$B^2_{(1)}$,
in the dual basis
are evaluated on
$\Kf^{(1)}$,
\begin{eqnarray}
    \Af_1^{(1)} 
&=&
    \left(
          \begin{array}{c}
          0\\
          1
          \end{array}
    \right)
    -
    \eps
    \left(
          \begin{array}{c}
          1\\
          \frac{\kappa (s+\kappa-\lambda)}{(s+\kappa)^3}
          \end{array}
    \right)
    +
    \eps^2
    \left(
          \begin{array}{c}
          0\\
          \frac{\kappa\lambda s (s+\kappa-\lambda)}
               {(s+\kappa)^6}
          \end{array}
    \right) 
    + 
    {\mathcal O}(\eps^3) , 
\Label{A1(1)}\\
    \Bf^2_{(1)} 
&=&
    \left(A_{21}^{(1)} , \, - A_{11}^{(1)}\right)
\Label{B2(1)}
\end{eqnarray}
Thus,
we see that
$A_1^{(1)}$
is tangent to the fast fibers 
at their base points
up to and including terms 
of ${\cal O}({\eps})$
as
Eq.~(\ref{MMH-TF1})
(with 
$\alpha= 1 $,
$\beta = -\frac{\kappa(s+\kappa-\lambda)}{(s+\kappa)^3}$)
implies.
As a result,
${\mathcal L}_\eps^{(1)}$
approximates ${\mathcal T}\Fe$
also up to and including terms 
of ${\cal O}({\eps})$.

\noindent\textbf{Remark 5.1.}
If,
in this particular example,
one evaluates
$A_1^{(1)}$ on $\Kf^{(0)}$
as opposed to $\Kf^{(1)}$
as we did above,
then 
the approximation of
${\mathcal T}\Fe$
is also accurate 
up to and including terms
of ${\mathcal O}(\eps)$.

\subsubsection{Second Iteration}
The blocks of 
$\bfL_{(1)}$
are
\begin{eqnarray}
    \bfL^{11}_{(1)}
&=&
    - (s+\kappa)
    + \eps \frac{(s+\kappa-\lambda)}{s+\kappa}
      \left[(c-1)
      + (c-\frac{s}{s+\kappa})\right]
    \nonumber\\
    &&\mbox{}
    + \eps^2 \frac{(c-1)(s+\kappa-\lambda)}{(s+\kappa)^3}
      \bigg[- \lambda (c-1)
       + [(s+\kappa-\lambda)c - s]\bigg] ,~~
\Label{L11 genl pt}\\
    \bfL^{12}_{(1)}
&=&
      \frac{s}{s+\kappa}-c
    + \eps \frac{c-1}{(s+\kappa)^2}
    \bigg[\lambda (c-1)
     - [(s+\kappa-\lambda)c - s]\bigg] ,
\Label{L12 genl pt}\\
    \bfL^{21}_{(1)}
&=&
    \left.\frac{\eps^2}{(s+\kappa)^2} 
    \right[(c-1)(s+\kappa-\lambda)(s+\kappa-2\lambda)
    \nonumber\\
    &&\left.\mbox{}
            + \lambda [(s+\kappa-\lambda)c - s]
            + (s+\kappa-\lambda)^2
              \left(
                    c-\frac{s}{s+\kappa}
              \right)\right] ,
\Label{L21 genl pt}\\
    \bfL^{22}_{(1)}
&=&
    \frac{\eps}{s+\kappa} 
    \left[\lambda (c-1)
          + (s+\kappa-\lambda)(\frac{s}{s+\kappa}-c)\right]
    \nonumber\\
    &&\mbox{}
    + \eps^2 \frac{(c-1)(s+\kappa-\lambda)}{(s+\kappa)^3}
    \bigg[\lambda (c-1) - [(s+\kappa-\lambda)c - s]\bigg] ,
\Label{L22 genl pt}
\end{eqnarray}
with remainders of $\mathcal{O}(\eps^3)$.

In the first step,
we update 
$A_2^{(1)}$
and
$B^1_{(1)}$
and evaluate 
the updated quantities on 
$\Kf^{(1)}$, 
to obtain
\begin{eqnarray}
%
    \Af_{12}^{(2)} 
&=&
      1
    + \eps^2 
    \frac{\kappa\lambda (2s-\kappa)(s+\kappa-\lambda)}{(s+\kappa)^6} ,
\Label{A12(2)}\\
    \Af_{22}^{(2)} 
&=&
      \frac{\kappa}{(s+\kappa)^2}
    + \eps \frac{\kappa\lambda(\kappa-3s)}{(s+\kappa)^5}
    \nonumber\\
    &&\mbox{}
    + \eps^2
    \frac{\kappa^2\lambda (7s-2\kappa)(s+\kappa-\lambda)+\kappa\lambda^2 s (s-2\kappa)}
    {(s+\kappa)^8} ,
\Label{A22(2)}\\
    \Bf^1_{(2)} 
&=&
    \left(- \Af_{22}^{(2)} , \, \Af_{12}^{(2)}\right) ,
\Label{B1(2)}
\end{eqnarray}
up to and including 
terms of $\mathcal{O}(\eps^2)$.

The CSP condition
\begin{eqnarray}
    \Bf^1_{(2)} g
&=&
      s
    - (s+\kappa)c
    - \eps\frac{\kappa(-s+(s+\kappa-\lambda)c)}{(s+\kappa)^2}
    \nonumber\\
    &&\mbox{}
    + \eps^2\kappa\lambda
    \left(  \frac{(3s-\kappa)(-s+(s+\kappa-\lambda)c)}
            {(s+\kappa)^5}\right.
    \nonumber\\
    &&\mbox{}
          \left.+ \frac{(2s-\kappa)(s+\kappa-\lambda)(s-(s+\kappa)c)}
            {(s+\kappa)^6}
    \right)
    + \mathcal{O}(\eps^3) \nonumber \\
&=& 
0 ,
\Label{f1(2)=0}
\end{eqnarray}
is satisfied if
\begin{eqnarray}
    c 
&=& 
    \frac{s}{s+\kappa}
+
    \eps\frac{\kappa\lambda s}{(s+\kappa)^4}
+    
    \eps^2
     \frac{\kappa\lambda s(2\kappa\lambda-3\lambda s-\kappa s-\kappa^2)}{(s+\kappa)^7}+\mathcal{O}(\eps^3). 
\Label{proj3}
\end{eqnarray}
Equation~(\ref{proj3})
defines $\Kf^{(2)}$,
the CSPM of order two,
which agrees with $\Me$
up to and including terms
of $\mathcal{O}(\eps^2)$;
recall Eq.~(\ref{Me-MMH}).

%
%
%
%
%

Then,
in the second step,
we update
$A_1^{(1)}$
and
$B^2_{(1)}$
to obtain
\begin{eqnarray}
    \Af_{11}^{(2)} 
&=&\mbox{}
    - \eps \frac{s+\kappa-\lambda}{s+\kappa}
    - \eps^2
      \frac{1}{(s+\kappa)^3}
      \left[
            (s+\kappa-\lambda)
            (s+\kappa-2\lambda)
            (c - 1)
      \right.
\nonumber\\
&&\mbox{}
      \left.
          + (s+\kappa-\lambda)^2
            \left(c - \frac{s}{s+\kappa}\right)
          + \lambda[(s+\kappa-\lambda)c - s]
      \right] .
\Label{A11(2)-genl}\\
    \Af_{21}^{(2)} 
&=&
      1
    + \eps \frac{(s+\kappa-\lambda)(c-1)}{(s+\kappa)^2} 
    + \eps^2
      \frac{1}{(s+\kappa)^4}
      \Bigg[
            (s+\kappa-\lambda)
            \Big[
            (s+\kappa-2\lambda)
            (c - 1)
\nonumber\\
&&\mbox{}
          + (s+\kappa-\lambda)
            \left(c - \frac{s}{s+\kappa}\right)
          + \lambda c
            \Big] 
          - \lambda s
      \Bigg] 
      \left(2c - \frac{2s+\kappa}{s+\kappa}\right) ,
\Label{A21(2)-genl}\\
    \Bf^2_{(2)} 
&=&
    \left(\Af_{21}^{(2)} , \, - \Af_{11}^{(2)}\right) ,
\Label{B2(2)-genl}
\end{eqnarray}
with remainders of ${\mathcal O}(\eps^3)$.
Evaluating 
these expressions
on $\Kf^{(2)}$,
we obtain
\begin{eqnarray}
    \Af_{11}^{(2)} 
&=&\mbox{}
    - \eps \frac{s+\kappa-\lambda}{s+\kappa}
    + \eps^2 
    \frac{\kappa(s+\kappa-2\lambda)(s+\kappa-\lambda)+\lambda^2 s}{(s+\kappa)^4} ,
\Label{A11(2)}\\
    \Af_{21}^{(2)} 
&=&
      1
    - \eps \frac{\kappa(s+\kappa-\lambda)}{(s+\kappa)^3} 
    \nonumber\\
    &&\mbox{}
    + \eps^2
    \frac{(s+\kappa-\lambda)(\kappa^2(s+\kappa-2\lambda)+\kappa\lambda s)+\kappa\lambda^2 s}
    {(s+\kappa)^6},
\Label{A21(2)}\\
    \Bf^2_{(2)} 
&=&
    \left(\Af_{21}^{(2)} , \, - \Af_{11}^{(2)}\right) ,
\Label{B2(2)}
\end{eqnarray}
with remainders of
$\mathcal{O}(\eps^3)$.
Therefore,
$\Af_1^{(2)}$ 
is tangent to
the fast fibers
at their base points
up to and including terms
of ${\cal O}(\eps^2)$,
according to
Eq.~(\ref{MMH-TF2})
(with 
$\alpha= 1 $,
$\beta = -\frac{\kappa(s+\kappa-\lambda)}{(s+\kappa)^3}$,
$\gamma
=
\frac{(s+\kappa-\lambda)(\kappa^2(s+\kappa-2\lambda)+\kappa\lambda s)+\kappa\lambda^2 s}
{(s+\kappa)^6}$),
and
${\mathcal L}_\eps^{(2)}$
is an ${\mathcal O}(\eps^2)$-accurate approximation to 
${\mathcal T}\Fe$.

\noindent\textbf{Remark 5.2.}
If one evaluates,
in this particular example,
$A_1^{(2)}$ on $\Kf^{(1)}$
instead of on $\Kf^{(2)}$
as we did above,
then 
the approximation of
${\mathcal T}\Fe$
is also accurate 
up to and including terms
of ${\mathcal O}(\eps^2)$.

\section{Linear Projection of Initial Conditions \label{s-disc}}
\setcounter{equation}{0}
The main result of
this article,
Theorem~\ref{t-CSPFq},
states that
after $q$ iterations
the CSP method
successfully 
identifies
${\cal T}\Fe$
up to and including
terms of ${\cal O}(\eps^{q+1})$,
where this approximation is
given explicitly by
$A_1^{(q)}$.
This information
is postprocessed to
project
the initial conditions
on the CSPM of
order $q$.
In this section,
we discuss
the accuracy 
and 
limitations of 
this linear projection.

Geometrically,
one knows from
Fenichel's theory
that any given
initial condition~$x_0$
sufficiently close to
$\Me$ lies on a
(generally nonlinear)
fiber $\Fe^p$
with base point $p$
on $\Me$.
Hence,
the ideal projection
would be
$\pi_F(x_0) = p$
(the subscript $F$
stands for fiber
or Fenichel)
and this is,
in general,
a nonlinear projection.

Within the framework of
an algorithm that
yields only
linearized information about
the fast fibers,
one must ask
how best to
approximate this ideal.
A consistent approach 
is to 
identify a point on
the slow manifold
such that
the approximate linearized fiber
through it
also goes through
the given initial condition.
This approach was used,
for example,
by Roberts~\cite{R-1999}
for systems with
asymptotically stable center manifolds,
where we note that
a different method is
first used to
approximate the center manifold.
Also,
this approach
is exact
in the special case that
the perturbed fast fibers are
hyperplanes which 
need not be
vertical.
In general,
if $x_0$ lies on
the linearized fiber
${\mathcal L}_\eps^{p_1}$
and if
$\pi_F(x_0) = p_2$,
then
the error $\|p_1 - p_2\|$ 
made by projecting linearly
is ${\cal O}(\eps)$ and
proportional to
the curvature of
the fiber
(see also \cite{R-1999}).

For fast--slow systems,
there is yet
another way to
linearly project
initial conditions on
the slow manifold.
One projects along
the approximate CSPF to
the space ${\cal T}_p\Fe$,
where $p$ is
the point on
the CSPM
that lies on
the same $\eps = 0$ fiber
as the initial condition.
This type of projection
is also consistent,
in the sense that
it yields an
exact result for
$\eps = 0$,
but has an error
of ${\cal O}(\eps)$
for $\eps > 0$.
Moreover,
it is
algorithmically simpler,
since
it does not involve
a search for
the base point of
the linearized fiber on
which the initial conditions lie.
However,
it has the disadvantage that
the projection
is not exact
in the special case that
the fast fibers are
(non-vertical) hyperplanes.


\appendix
\section{The Action of the ${\cal O}(1)$ Jacobian
on ${\cal T}_p\Mo$ \label{s-lemma}}
\setcounter{equation}{0}
The spaces
${\cal T}_p\Fe$
and
${\cal T}_p\Me$
depend,
in general,
on both
the point $p\in\Me$
and $\eps$.
As a result,
the basis
$A$
also depends 
on $p$ and $\eps$,
and hence
$A_f$ and $A_s$
possess
formal asymptotic expansions
in terms of
$\eps$,
\begin{eqnarray}
    A_f = \sum_{i=0}\eps^i A^i_f ,
\quad
    A_s = \sum_{i=0}\eps^i A^i_s .
\Label{A-asympt-exp}
\end{eqnarray}
Next,
we compute
the action of
the Jacobian on
$A_s$
to leading order.

\begin{LEMMA}
\label{l-DgoAo}
Ker$(Dg(p))_0
=
{\cal T}_p\Mo$,
for $p\in\Mo$.
In particular,
$(Dg)_0 A^0_s = 0$.
\end{LEMMA}

\begin{PROOF}
The Jacobian is
a linear operator,
so 
it suffices to
show that
every column vector of
a basis for ${\cal T}_p\Mo$
vanishes under
the left action of
the Jacobian.
We choose
this basis to be
the matrix
$\left(\begin{array}{c}I_m\\D_yh_0\end{array}\right)$.

We compute
\begin{eqnarray}
    Dg_0 \left(
             \begin{array}{c}
             I_m\\
             D_yh_0
             \end{array}
           \right)
=
    \left(
          \begin{array}{cc}
          0 & 0\\
          D_yg_2 & D_zg_2
          \end{array}
    \right)
    \left(
          \begin{array}{c}
          I_m\\
          D_yh_0
          \end{array}
    \right)
=
    \left(
          \begin{array}{c}
          0\\
          D_yg_2 + D_zg_2 D_yh_0
          \end{array}
    \right) .~~
\Label{DgoTMo-aux}
\end{eqnarray}
Differentiating
both members of
the ${\cal O}(1)$
invariance equation
$g_2(y,h_0(y),0)=0$
with respect to
$y$,
we obtain
\be
    D_yg_2(y,h_0(y),0)
    +
    D_zg_2(y,h_0(y),0) D_yh_0(y)
    =
    0 .
\Label{0inv-id}
\ee
Equations~(\ref{DgoTMo-aux})
and (\ref{0inv-id})
yield the desired result
\begin{eqnarray}
    Dg_0 \left(
             \begin{array}{c}
             I_m\\
             D_yh_0
             \end{array}
           \right)
    &=& 
    \left(
             \begin{array}{c}
             0\\
             0
             \end{array}
           \right) ,
\hskip0.1truein\mbox{on $\Mo$} .
\Label{DgoTMo}
\end{eqnarray}

Finally,
the identity
$(Dg)_0A_s^0 = 0$
follows from
the fact that
$A^0_s$ spans
${\cal T}_p\Mo$,
since 
$A^0_s
=
A_s\vert_{\eps=0}$
by Eq.~(\ref{A-asympt-exp}).
\end{PROOF}


\begin{center}
\textbf{ACKNOWLEDGEMENTS}
\end{center}
\noindent
The work of H.~K.\ was supported by the Mathematical, Information, and
Computational Sciences Division subprogram of the Office of Advanced
Scientific Computing Research, U.S.~Department of Energy, under
Contract W-31-109-Eng-38.
The work of T.~K.\ and A.~Z.\ was supported in part 
by the Division of Mathematical Sciences
of the National Science Foundation via grant NSF-0306523.

\newpage
\noindent
Corresponding author:

\noindent
Hans G.\ Kaper \\
Division of Mathematical Sciences \\
National Science Foundation \\
4201 Wilson Boulevard, Suite 1025 \\
Arlington, VA 22230

\noindent
Authors' e-mail addresses:

\noindent
\texttt{azagaris@math.bu.edu} \\
\texttt{hkaper@nsf.gov, kaper@mcs.anl.gov} \\
\texttt{tasso@math.bu.edu} \\

\vfill

\begin{flushright}
\scriptsize
\framebox{\parbox{2.4in}{
The submitted manuscript has been created
by the University of Chicago as Operator of
Argonne National Laboratory ("Argonne")
under Contract No.\ W-31-109-ENG-38
with the U.S.\ Department of Energy.
The U.S.\ Government retains for itself,
and others acting on its behalf,
a paid-up, nonexclusive, irrevocable
worldwide license in said article
to reproduce, prepare derivative works,
distribute copies to the public,
and perform publicly and display publicly,
by or on behalf of the Government.}}
\normalsize
\end{flushright}


\begin{thebibliography}{99}
%
\bibitem{DFN-1985}
B.~A.~Dubrovin, A.~T.~Fomenko, and S.~P.~Novikov, 
\textsl{Modern Geometry -- Methods and Applications, Vol.~2},
Graduate Texts in Mathematics, \textbf{104},
Springer-Verlag, New York, 1985
%
\bibitem{F-1979}
N.~Fenichel,
Geometric singular perturbation theory
for ordinary differential equations,
{\it J.~Diff.\ Eq.}\ \textbf{31} (1979) 53--98
%
\bibitem{GL-1992}
D.~A.~Goussis and S.~H.~Lam,
A study of homogeneous methanol
oxidation kinetics using CSP,
in:
\textsl{Twenty-Fourth Symposium (International) on Combustion,
The University of Sydney,
Sydney, Australia, July 5--10, 1992},
The Combustion Institute, Pittsburgh, 1992, pp.~113--120
%
\bibitem{HG-1999}
M.~Hadjinicolaou and D.~A.~Goussis,
Asymptotic solutions of stiff PDEs with the CSP method:
The reaction diffusion equation,
{\it SIAM J.~Sci.\ Comput.}\ \textbf{20} (1999) 781--810
%
\bibitem{J-1994}
C.~K.~R.~T.~Jones,
Geometric singular perturbation theory,
in:
\textsl{Dynamical Systems, Montecatini Terme,} L.~Arnold,
Lecture Notes in Mathematics, \textbf{1609},
Springer-Verlag, Berlin, 1994, pp.~44--118
%
\bibitem{L-1993}
S.~H.~Lam,
Using CSP to understand complex chemical kinetics,
{\it Combust.\ Sci.\ Tech.}\ \textbf{89} (1993) 375--404
%
\bibitem{LG-1988}
S.~H.~Lam and D.~A.~Goussis,
Understanding complex chemical kinetics with computational
singular perturbation,
in {\it Twenty-Second Symposium (International) on Combustion,
The University of Washington, Seattle, Washington, August 14--19,
1988,} The Combustion Institute, Pittsburgh, 1988, pp. 931--941
%
\bibitem{LG-1991}
S.~H.~Lam and D.~A.~Goussis,
Conventional asymptotics 
and computational singular perturbation theory
for simplified kinetics modeling,
in {\it Reduced Kinetic Mechanisms 
and Asymptotic Approximations
for Methane-Air Flames}, 
M.~Smooke, ed.,
Lecture Notes in Physics \textbf{384},
Springer-Verlag, New York, 1991,
Chapter 10
%
\bibitem{LG-1994}
S.~H.~Lam and D.~A.~Goussis,
The CSP method for simplifying kinetics,
{\it Internat. J.~Chem.\ Kin.}\ \textbf{26} (1994) 461--486
%
%
\bibitem{LJL-2001}
T.~F.~Lu, Y.~G.~Ju, and C.~K.~Law,
Complex CSP for chemistry reduction and analysis,
{\it Combustion and Flame}\ \textbf{126} (2001) 1445--1455
%
\bibitem{MDMG-1999}
A.~Massias, D.~Diamantis, E.~Mastorakos, and D.~Goussis,
Global reduced mechanisms for methane and hydrogen combustion
with nitric oxide formation constructed with CSP data,
{\it Combust.\ Theory Modelling} \textbf{3} (1999) 233--257
%
\bibitem{MG-2001}
A.~Massias and D.~A.~Goussis,
On the manifold
of stiff
reaction-diffusion PDE's:
The effects of diffusion,
preprint (2001)
%
\bibitem{M-1996}
K.~D.~Mease,
Geometry of computational singular perturbations,
in \textsl{Nonlinear Control System Design,}
vol. 2,
A.~J.~Kerner and D.~Q.~M.~Mayne, editors,
Pergamon Press, Oxford,  U.K., 1996, pp.~855--861
%
\bibitem{O-1986}   
P.~J.~Olver,
\textsl{Applications of Lie Groups to Differential Equations},
Graduate Texts in Mathematics, \textbf{107},
Springer-Verlag, New York, 1986
%
\bibitem{P-1987}
B.~O.~Palsson,
On the dynamics of the irreversible
Michaelis--Menten reaction mechanism,
{\it Chem.\ Eng.\ Sci.}\ \textbf{42} (1987) 447--458
%
\bibitem{PL-1984}
B.~O.~Palsson and E.~N.~Lightfoot,
Mathematical modelling of dynamics and control
in metabolic networks. I.~On Michaelis--Menten kinetics,
{\it J.~theor.\ Bio.}\ \textbf{111} (1984) 273--302
%
\bibitem{R-1999}
A.~J.~Roberts,
Computer algebra derives
correct initial conditions
for low-dimensional
dynamical systems,
arXiv: chao-dyn/9901010
%
\bibitem{VG-2001}
M.~Valorani and D.~A.~Goussis,
Explicit time-scale splitting algorithm
for stiff problems: auto-ignition of
gaseous mixtures behind a steady shock,
{\it J.~Comp.\ Phys.}\ \textbf{169} (2001) 44--79
%
\bibitem{VNG-2003}
M.~Valorani, H.~M.~Najm, and D.~A.~Goussis,
CSP analysis of 
a transient 
flame-vortex interaction: 
time scales and manifolds,
{\it Combustion and Flame}\ \textbf{134} (2003) 35--53
%
\bibitem{ZKK-2003}
A.~Zagaris, H.~G.~Kaper, and T.~J.~Kaper,
Analysis of the 
Computational Singular Perturbation
reduction method for
chemical kinetics,
{\it J. Nonlin. Sci.}
(to appear);
also available at
arXiv: math.DS/0305355
\end{thebibliography}
\end{document}